\documentclass[a4paper,11pt]{article}
\usepackage[margin=0.6in]{geometry}
\pdfoutput=1
\usepackage{mathtools}

\usepackage[utf8]{inputenc}
\usepackage{amsfonts,amssymb, bm,amsthm}
\usepackage{graphicx,color,epstopdf}
\usepackage[inline]{showlabels}
\allowdisplaybreaks
\newtheorem{mylemma}{Lemma}[section]
\newtheorem{mytheorem}{Theorem}[section]

\newtheorem{myremark}{Remark}[section]

\def\XXint#1#2#3{{\setbox0=\hbox{$#1{#2#3}{\int}$}
    \vcenter{\hbox{$#2#3$}}\kern-.5\wd0}}

\def\bi{{\bf i}}

\begin{document}
\title{A Fourier-matching Method for Analyzing Resonance
  Frequencies by a Sound-hard Slab with Arbitrarily Shaped Subwavelength Holes}
\author{Wangtao Lu$^1$, Wei Wang$^2$ and Jiaxin Zhou$^3$}
\footnotetext[1]{School of Mathematical Sciences, Zhejiang University, Hangzhou
  310027, China. Email: wangtaolu@zju.edu.cn. This author
  is partially supported by NSF of Zhejiang Province for Distinguished Young
  Scholars (LR21A010001).} \footnotetext[2]{School of Mathematical Sciences,
  Zhejiang University, Hangzhou 310027, China. Email: wangw07@zju.edu.cn.}
\footnotetext[3]{School of Mathematical Sciences, Zhejiang University, Hangzhou
  310027, China. Email: jiaxinzhou@zju.edu.cn.}
\maketitle
\begin{abstract}
  This paper presents a simple Fourier-matching method to rigorously study
  resonance frequencies of a sound-hard slab with a finite number of arbitrarily
  shaped cylindrical holes of diameter ${\cal O}(h)$ for $h\ll1$. Outside the
  holes, a sound field can be expressed in terms of its normal derivatives on
  the apertures of holes. Inside each hole, since the vertical variable can be
  separated, the field can be expressed in terms of a countable set of Fourier
  basis functions. Matching the field on each aperture yields a linear system of
  countable equations in terms of a countable set of unknown Fourier
  coefficients. The linear system can be reduced to a finite-dimensional linear
  system based on the invertibility of its principal submatrix, which is proved
  by the well-posedness of a closely related boundary value problem for each hole
  in the limiting case $h\to 0$, so that only the leading Fourier coefficient of
  each hole is preserved in the finite-dimensional system. The resonance
  frequencies are those making the resulting finite-dimensional linear system
  rank deficient. By regular asymptotic analysis for $h\ll1$, we get a systematic
  asymptotic formula for characterizing the resonance frequencies by the 3D
  subwavelength structure. The formula reveals an important fact that when all
  holes are of the same shape, the $Q$-factor for any resonance frequency
  asymptotically behaves as ${\cal O}(h^{-2})$ for $h\ll1$ with its prefactor
  independent of shapes of holes.
\end{abstract}
\section{Introduction}
Subwavelength structures have attracted great attentions in the area of wave
scattering problems in the past decades \cite{astlalpal00,carmahgar11, chenet13, ebblezwol98, genebb07, garmarebbkui10, liulal08, seoetal09, stupodgor10, tak01, yansam02}. These structures have been
experimentally observed and numerically simulated to own some exclusive
features, such as extraordinary optical transmission, local field enhancement,
making themselves widely applicable in areas such as biological sensing and
imaging, microscopy, spectroscopy and communication \cite{sarvig07,lienyllud83}.
It has now been well-known that these features are mostly caused by the
existence of high-Q resonances in subwavelength structures. Mathematically, a
resonance frequency $k$ can be defined as a complex frequency in the lower-half
of complex plane $\mathbb{C}$, at which the scattering problem loses uniqueness.
The quality factor defined as $Q={\rm Re}(k)/(2{\rm Im}(k))$ can be used to
measure how great wave field can be enhanced in subwavelength structures at the
real frequency ${\rm Re}(k)$. Therefore, it is highly desired to design a
subwavelength structure with a resonance frequency closing enough to the real
axis.

To this purpose, existing literatures have made great efforts in the past to
propose either effectively computational methods or rigorously mathematical
theories to quantitatively analyze resonance frequencies in subwavelength
structures \cite{babbontri10,bonsta94,bontri10,
  braholsch20,cladurjoltor06,gaoliyua17,holsch19,huyualu20,joltor06a,joltor06b,joltor08,linshizha20,linzha17,linzha18a,linzha18b,
  shivol07,shi10}. Among existing theories, roughly two types of methods have
been proposed: boundary-integral-equation (BIE) method and matched-asymptotics
method, mainly to study two-dimensional (2D) subwavelength structures. Bonnetier
and Triki \cite{bontri10} used the first method to firstly study wave scattering
by a perfectly conducting half plane with a subwavelength cavity and obtained an
asymptotic formula of resonance frequencies. Subsequently, Babadjian et al.
\cite{babbontri10} used this method to study resonances by two interacting
subwavelength cavities; Lin and Zhang developed a simplified BIE method to study
resonances by a slab with a single 2D slit \cite{linzha17}, periodic slits
\cite{linzha18a,linzha18b}, or a periodic array of two subwavelength slits
\cite{linshizha20}; Gao et al. \cite{gaoliyua17} studied resonance frequencies
by a rectangular cavity with different conducting boundaries. Using the second
method, Joly and Tordeux \cite{joltor06a,joltor06b,joltor08} and Clausel et al.
\cite{cladurjoltor06} studied resonances by thin slots; Holley and Schnitzer
\cite{holsch19} studied resonances of a slab with a single slit, and Brand\~ao
et al. \cite{braholsch20} studied resonances of a slab of finite conductivity
with a single slit or periodic slits.

Compared with 2D structures, three-dimensional (3D) subwavelength structures are
more flexible in practical fabrication and in fact are much easier to realize
high-Q resonators \cite{genebb07, garmarebbkui10, liulal08, carmahgar11}.
Nevertheless, much fewer theories have been developed so far to rigorously study
resonances of 3D structures \cite{liazou20}. In a recent work \cite{zhlu21}, the
authors developed a Fourier-matching method to study resonances by a slab of
finite number of 2D slits. Unlike existing methods, the Fourier-matching method
does not use the complicated Green function in a slit so that the overall theory
becomes more straightforward. Consequently, we are motivatied to extend this
simple Fourier-matching method to 3D subwavelength structures. Inheriting its
advantages, this paper further simplifies the original analyzing procedure of
the Fourier-matching method, making it applicable for studying resonances of a
sound-hard slab with a finite number of 3D subwavelength cylindrical holes of
arbitrary shapes, as discussed below.

\begin{figure}[!ht]
  \centering
  a)\includegraphics[width=0.42\textwidth]{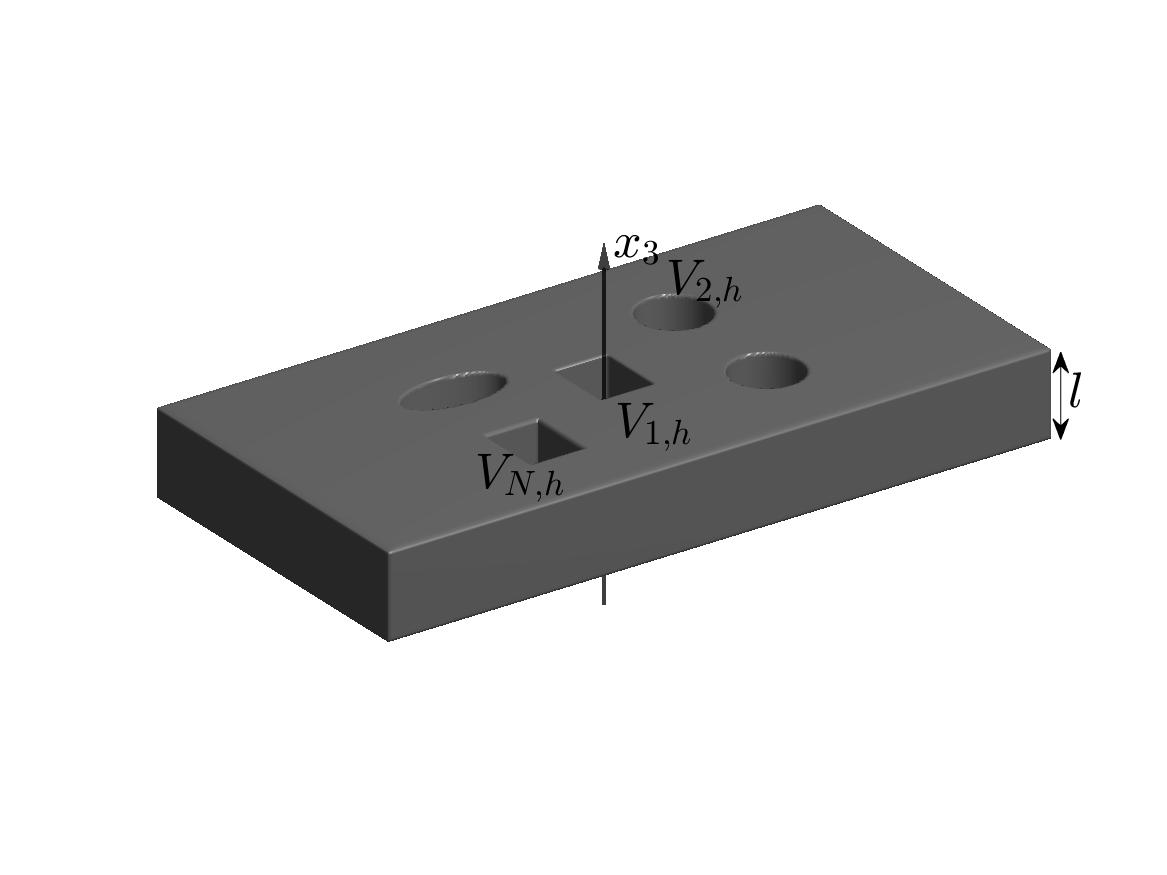}
  b)\includegraphics[width=0.4\textwidth]{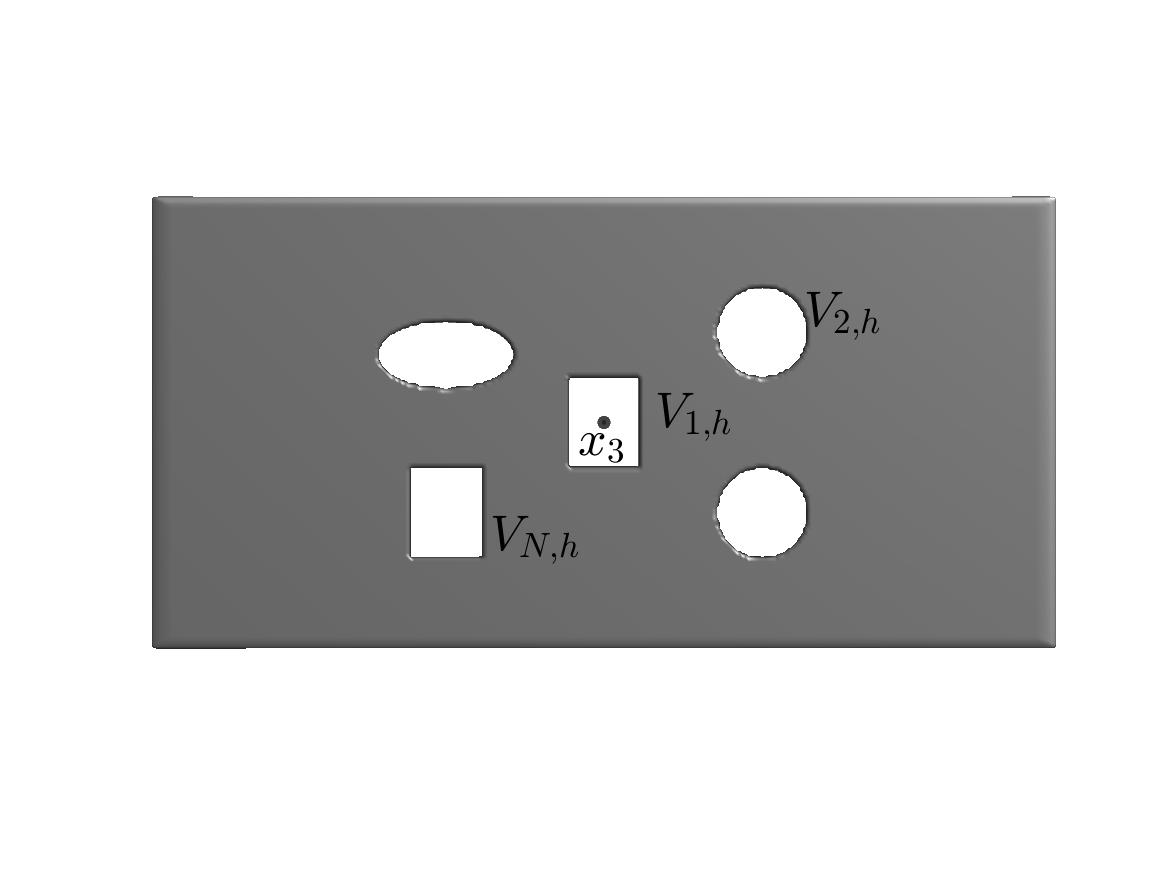}
  %\vspace{-1.5cm}
  \caption{A sound-hard slab of thickness $l$ with $N$ cylindrical holes $\{V_{j,h}\}_{j=1}^N$: (a) side
    view; (b) top-view.}
  \label{fig:model}
\end{figure}

As shown in Figure~\ref{fig:model},
let $\{V_{j,h}\}_{j=1}^N$ denotes the set of $N$ holes in a sound-hard slab of
thickness $l$. Throughout this paper, we assume that $\{V_{j,h}\}_{j=1}^N$
satisfies the following conditions: (1) $\{V_{j,h}\}_{j=1}^N$ are cylindrical,
i.e. $x_3$-independent; (2) $\{V_{j,h}\}_{j=1}^N$ are generated respectively by
two-dimensional, simple-connected Lipschitz domains $\{G_j\}_{j=1}^N$, all of
which contain the origin point $(0,0)$, and the small parameter $h\ll 1$,
through
\[
  V_{j,h} = \{(x_1,x_2,x_3)\in \mathbb{R}^3:(x_1,x_2)\in D_j+hG_j,x_3\in(-l,0)\},
\]
where $D_j+hG_j=\{D_j+(hx_1,hx_2)\in\mathbb{R}^2:(x_1,x_2)\in G_j\}$,
$\{D_j\}_{j=1}^{N}$ are well-separated points in $\mathbb{R}^2$ so that
$C_j=(D_j,-l/2)$ become the center of $V_{j,h}$; (3) the area of $G_j$ is
$1$. Here, condition (3) is not necessary and is only introduced to simplify the
presentation. Let $\Gamma_{j,h}=\partial V_{j,h}\cap\{x:x_3=0\}$ denotes the top
boundary of $V_{j,h}$, which will be called the aperture of $V_{j,h}$ in the following.

In such a 3D structure, a sound field outside the holes $\{V_{j,h}\}$ can be
expressed in terms of its normal derivatives on the apertures $\{G_{j,h}\}$.
Inside each hole $V_{j,h}$, since the vertical variable can be separated, the
field can be expressed in terms of a countable set of Fourier basis functions.
Matching the field on each aperture $G_{j,h}$ yields a linear system of
countable equations in terms of a countable set of unknown Fourier coefficients.
The linear system can be reduced to a finite-dimensional linear system based on
the invertibility of a principal submatrix, which is proved by establishing the
well-posedness of a closely related boundary value problem for each hole
$V_{j,h}$ in the limiting case $h\to 0$, so that only the leading Fourier
coefficient for each hole $V_{j,h}$ is preserved in the finite-dimensional
system. In other words, the resonance frequencies are those making the resulting
$N$-dimensional linear system rank deficient. By regular asymptotic analysis, we
get a systematic asymptotic formula for characterizing the resonance frequencies
of the 3D structure. The formula reveals an important fact that when all holes
are of the same shape, the quality factor $Q$ for any resonance frequency
asymptotically behaves as ${\cal O}(h^{-2})$ for $h\ll1$, the prefactor of which
depends only on the locations of holes, but is independent of the shape.

\subsection{Notations and Equivalent Sobolev norms}

For any Lipschitz domain $\Omega\in \mathbb{R}^n, n=2,3$, let $L^2(\Omega)$ denote the set of all
square-integrable functions on $\Omega$ equipped with the usual inner product:
for any $f,g\in L^2(\Omega)$,
\[
(f,g)_{L^2(\Omega)}=\int_{\Omega}f\bar{g}dx.
\]
Let $H^1(\Omega)=\{f\in L^2(\Omega): \nabla f\in
(L^2(\Omega))^3\}$ be equipped with the standard inner product: for any $f,g\in
H^1(\Omega)$,
\[
  (f,g)_{H^1(\Omega)} = \int_{\Omega}f\bar{g} + \nabla f\cdot\nabla\bar{g}dx.
\]
Let $\tilde{H}^{-1}(\Omega)$ denote the dual space of $H^1(\Omega)$. Then, the
interpolation theory can help us to define the fractional Sobolev space
$H^{1/2}(\Omega)$ and its dual space $\tilde{H}^{-1/2}(\Omega)$. Let
$\ell^2=\{\{f_{m}\}_{m=0}^{\infty}\subset\mathbb{C}:
\sum_{m=0}^{\infty}|f_{m}|^2<\infty\}$ be equipped with its natural inner
product: for any $\{f_m\}_{m=0}^{\infty},\{g_m\}_{m=0}^{\infty}\in\ell^2$:
\[
  <\{f_m\},\{g_m\}>_{\ell^2} = \sum_{m=0}^{\infty}f_m\bar{g}_m.
\]

For any $f\in H^{1/2}(\Omega)$ and $g\in \tilde{H}^{-1/2}(\Omega)$, the duality
pair $<f,g>_{H^{1/2}(\Omega),\tilde{H}^{-1/2}(\Omega)}$, understood as the
functional $f$ acting on $g$, will be abbreviated as $<f,g>_{1/2,-1/2}$ for
simplicity when the definition domain of $f$ or $g$ is clear from the context;
similarly, $<g,f>_{-1/2,1/2}$ denotes the functional $g\in
\tilde{H}^{-1/2}(\Omega)=( H^{1/2}(\Omega) )'$ acting on element $f\in
H^{1/2}(\Omega)$. Certainly, $<f,g>_{1/2,-1/2}=\overline{<g,f>_{-1/2,1/2}}$ and
becomes $(f,g)_{L^2(\Omega)}$ when $g\in L^2(\Omega)$.

To simply characterize the aforementioned fractional Sobolev spaces on the
aperture boundary $\Gamma_{j,h}$ of the hole $V_{j,h}$, we rely on the following
theorem regarding spectral properties of the 2D Laplacian
$\Delta_2=\partial_{x_1}^2+\partial_{x_2}^2$ on the Lipschitz domains
$\{G_j\}_{j=1}^N$.
\begin{mytheorem}[See Theorem 4.12 in \cite{mcl00}]
  \label{thm:basis}
  For the 2D Lipschitz domain $\{G_j\}_{j=1}^N$ generating the $N$
  holes $\{V_{j,h}\}_{j=1}^N$, respectively, there exist sequences
  of functions $\phi_{1,j},\phi_{2,j},\cdots,$ in $H^1(G_j)$, and of nonnegative
  numbers $\lambda_{0,j}, \lambda_{1,j},\cdots,$ having the following properties:
  \begin{itemize}
  \item[(i)] Each $\phi_{m,j}$  is an eigenfunction of $-\Delta_2$ with eigenvalue $\lambda_{m,j}$:
    \begin{align*}
      -\Delta_2\phi_{m,j} =& \lambda_{m,j}\phi_{m,j},\quad{\rm on}\quad G_{j},\\
      \partial_{\nu} \phi_{m,j} =& 0,\quad{\rm on}\quad\partial G_{j}.
    \end{align*}
  \item[(ii)] The eigenvalues satisfy $0=\lambda_{0,j}<\lambda_{1,j}\leq
    \lambda_{2,j}\leq \cdots$ with $\lambda_{m,j}\to\infty$ as $i\to\infty$.
  \item[(iii)] The eigenfunctions $\{\phi_{m,j}\}_{i=0}^{\infty}$ form a complete
    orthonormal system in $L_2(G_j)$, and in particular $\phi_{0,j}=1$.
  \item[(iv)] For any $f\in H^1(G_j)$,
    \[
      ||f||_{H^1(G_j)}^2\eqsim \sum_{m=0}^{\infty}(1+\lambda_{m,j})|(f,\phi_{m,j})_{L^2(G_j)}|^2.
    \]
  \end{itemize}
  \begin{proof}
    The only difference from Theorem 4.12 of \cite{mcl00} is $\lambda_{0,j}=0$,
    $\lambda_{1,j}>0$ and $\phi_{0,j}=1$. This can be seen that all eigenvalues
    must be nonnegative by testing the governing equation with $\phi_{m,j}$
    themselves. On the other hand, for $\lambda_{0,j}=0$, $\phi_{0,j}=1$ is the
    unique up to a sign and normalized (see condition (3) about $G_j$) solution
    of the eigenvalue problem so that $\lambda_{m,j}>0$ for $i\geq 1$.
  \end{proof}
\end{mytheorem}

Based on Theorem~\ref{thm:basis}(iii), on each aperture $\Gamma_{j,h}=(D_j+hG_{j})\times\{x_3=l\}$,
\[
  \{\phi_{m,j}(\cdot;h)=h^{-1}\phi_{m,j}((\cdot-D_j)/h)\}_{m=0}^{\infty},
\]
forms a complete and orthonormal basis of $L^2(\Gamma_{j,h})$, so that for any
$f\in L^2(\Gamma_{j,h})$, the set of Fourier coefficients,
\[\{f_{m,j} =
  (f,\phi_{m,j}(\cdot;h))_{L^2(\Gamma_j,h)}\}_{m=0}^{\infty}\in\ell^2.\]
By Parserval's identity,
\[
  ||f||_{L^2(\Gamma_{j,h})} = ||\{f_m\}||_{\ell^2}=\left(  \sum_{m=0}^{\infty}|f_m|^2\right)^{1/2}.
\]
On the other hand, Theorem~\ref{thm:basis}(iv) indicates that,  the $H^1(\Gamma_{j,h})$ can be equipped
with the following equivalent norm: for any $f\in H^1(\Gamma_{j,h})$,
\[
  ||f||_{H^1(\Gamma_{j,h})} = \left(  \sum_{m=0}^{\infty}(1+\lambda_m)|f_m|^2\right)^{1/2}.
\]
Now by interpolation theory, $H^{1/2}(\Gamma_{j,h})$ can be equipped with the
following equivalent norm
\[
  ||f||_{H^{1/2}(\Gamma_{j,h})} = \left(\sum_{m=0}^{\infty}(1+\lambda_{m,j})^{1/2}|f_m|^2 \right)^{1/2},
\]
so that its dual space $\tilde{H}^{-1/2}(\Gamma_{j,h})$ is equipped with
\[
  ||f||_{\tilde{H}^{-1/2}(\Gamma_{j,h})} = \left(  \sum_{m=0}^{\infty}(1+\lambda_{m,j})^{-1/2}|f_m|^2\right)^{1/2},
\]
where $f_m$ should be redefined as $<f,\phi_m>_{-1/2,1/2}$ now.

The rest of this paper is organized as follows. In section 2, we study
resonances by a sound-soft slab with a single hole, and analyze the mechanism of
field enhancement near resonance frequencies. In section 3, we extend the
approach to study resonances of a slab with multiple slits and provide an
accurate asymptotic formula of the resonance frequencies. Finally, we draw the
conclusion and present some potential applications of the current method.

\section{Single cylindrical hole}
To clarify the basic idea, we begin with a slab of thickness $l$ with a single
cylindrical hole, say $V_{1,h}$, for $h\ll1$. For
simplicity, we assume $D_1=(0,0)$ so that $C_1=(0,0,-l/2)$ in this section. We
seek a complex frequency $k$ such that there exists a nonzero outgoing sound
wave field $u$ satisfies the following three-dimensional Helmholtz equation
\begin{align}
  \Delta u + k^2 u &= 0,\quad {\rm on}\ \Omega_h,\\
  \partial_{\nu} u &= 0,\quad {\rm on}\ \partial\Omega_h,
\end{align}
where $\Omega_h$ is the interior of $\{x\in\mathbb{R}^3: x_1,x_2\in \mathbb{R},
x_3\notin[-l,0]\}\cup \overline{V_{1,h}}$. The searching region ${\cal
  B}=\{k\in\mathbb{C}: {\rm Re}(k)>0, {\rm Im}(k)<0, |k|\in (\epsilon_0,K)\}$
for some sufficiently small constant $\epsilon_0>0$ and some sufficiently large
constant $K>0$. Similar as in \cite{zhlu21}, we consider even modes and odd modes
symmetric about $x_3=-l/2$.
\subsection{Even modes}
\begin{figure}[!ht]
  \centering
  a)\includegraphics[width=0.35\textwidth]{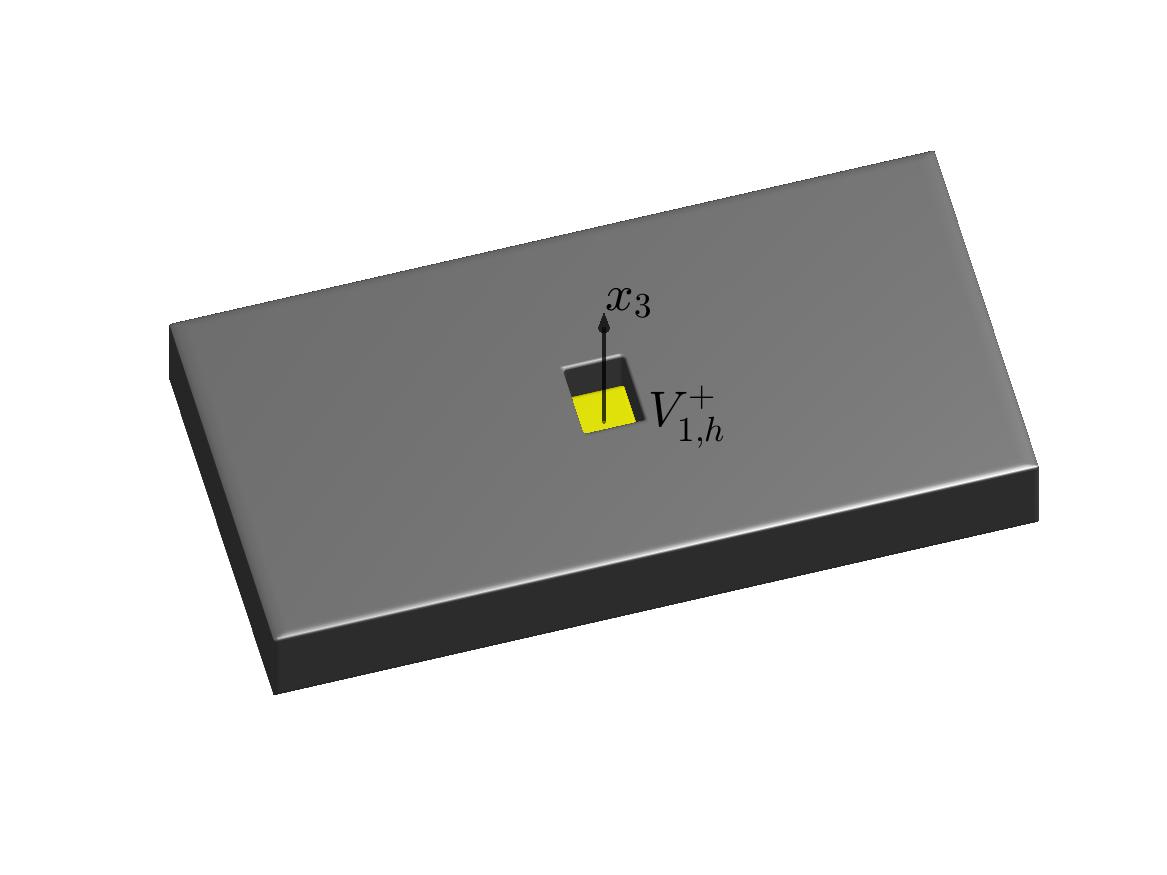}
  b)\includegraphics[width=0.4\textwidth]{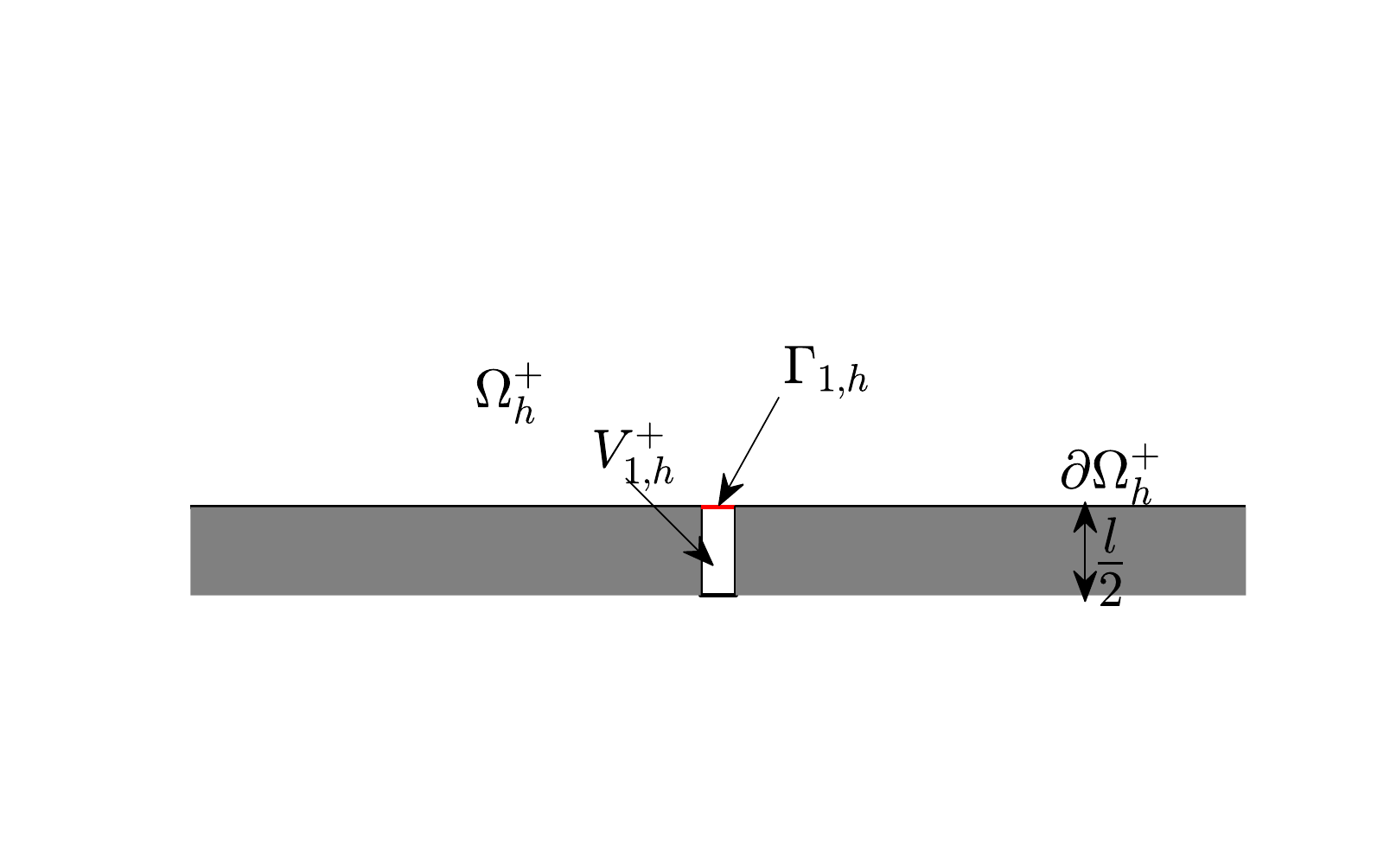}
  %\vspace{-1.5cm}
  \caption{A half-plane problem with one cylindrical cavity $V_{1,h}^+$: (a) side
    view; (b) the cross-section at $x_1=0$.}
  \label{fig:model2}
\end{figure}

Suppose now $u$ is an even function about $x_3=-l/2$, i.e., $u(x_1,x_2,x_3)=
u(x_1,x_2,-l-x_3)$, so that the original
problem is reduced to the following half-space problem: find a solution $u\in
H^{1}_{\rm loc}(\Omega_h^+)$ solving
\begin{align}
  \label{eq:gov:e}
  \Delta u + k^2 u &= 0,\quad {\rm on}\ \Omega_h^+,\\
  \label{eq:bc:e}
  \partial_{\nu} u &= 0,\quad {\rm on}\ \partial\Omega_h^+,
\end{align}
where $\Omega_h^+=\mathbb{R}_+^3\cup V_{1,h}^+\cup \Gamma_{1,h}$,
$V_{1,h}^+=V_{1,h}\cap\{x\in\mathbb{R}^3:x_3\in(-l/2,0)\}$ and we recall that
$\Gamma_{1,h}$ denotes the aperture surface, as shown in
Figure~\ref{fig:model2}.
Recall from Theorem~\ref{thm:basis} that
\[
  \{\phi_{m,1}(x;h)=h^{-1}\phi_{m,1}((x-C_1)/h)\}_{m=0}^{\infty},
\]
form a complete and orthonormal basis in $L^2(\Gamma_h)$, and the corresponding
eigenvalues are $\{\lambda_{m,1}\}_{m=2}^{\infty}$. To simplify
the presentation, we shall suppress the label $1$ so that $V_h^+=V_{1,h}^+$,
etc..

Since $u|_{\Gamma_h}\in H^{1/2}(\Gamma_h)$, there exists a sequence $\{b_{m}\}_{m=0}^{\infty}$ such that
\begin{align}
  \label{eq:rep:u}
  u|_{\Gamma_h} = \sum_{m=0}^{\infty}b_{m}\phi_{m}(x;h)[e^{\bi s_{m} l} + 1],
  %u|_{\Gamma_A} = \sum_{m,n=0}^{\infty}b_{mn}\phi_m(x_1)\phi_n(x_2)[e^{\bi s_{mn}(x_3+l)} + e^{-\bi s_{mn} x_3}]
\end{align}
with
\[
  \sum_{m=0}^{\infty}(1+\lambda_m)^{1/2}|b_{m}(e^{\bi s_{m} l}+1)|^2<\infty,
\]
where
\begin{align}
  \label{eq:sm}
  s_{m}=& \sqrt{k^2 - \frac{\lambda_m}{h^2}}.
\end{align}
Therefore, $\{a_{m}=\lambda_m^{1/2}b_{m}\}_{m=0}^{\infty}\in \ell^2$ for $h\ll 1$ since
\[
  \sum_{m=0}|a_{m}|^2 \leq C\sum_{m=0}^{\infty}(1+\lambda_m)^{1/2}|b_{m}(e^{\bi s_{m} l}+1)|^2<\infty.
\]
In $V_{h}^+$, define
\begin{equation}
  \label{eq:u-}
  u^-(x) = \sum_{m=0}^{+\infty}b_{m}\phi_{m}(x;h)[e^{\bi s_{m}(x_3+l)} + e^{-\bi s_{m} x_3}].
\end{equation}
We could verify that $\phi=u-u^-\in H^1(V_h^+)$ solves
\begin{align}
  \Delta \phi =& -k^2 \phi,\quad{\rm on}\quad V_h^+,\\
  \partial_{\nu }\phi =& 0,\quad{\rm on}\quad \partial V_h^+\backslash\overline{\Gamma_h},\\
  \phi=&0,\quad{\rm on}\quad\Gamma_h.
\end{align}
When $k\in {\cal B}$ and $h\ll 1$, $-k^2$ is not an eigenvalue of the above
problem so that we have $u=u^-$ on $C_h^+$. Therefore, its normal derivative on
$\Gamma_h$ becomes
\begin{equation}
  \label{eq:normder}
  \partial_{\nu}u(x) = \partial_{x_3}u(x_1,x_2,0) = \sum_{m=0}^{+\infty}b_{m}\bi
  s_{m}\phi_{m}(x;h) [e^{\bi s_{m}l} - 1]\in \tilde{H}^{-1/2}(\Gamma_h).
\end{equation}
When $e^{\bi s_{0}l}+1 = 0$, the representation (\ref{eq:rep:u}) of
$u|_{\Gamma_h}$ is impossible to define $b_{0}$. To resolve this issue, we
could choose (\ref{eq:normder}) to define all $\{b_{m}\}_{m=0}^{\infty}$ so that
the representation (\ref{eq:u-}) becomes valid for all finite frequencies
$k\in{\cal B}$ and $h\ll 1$.

Let $\epsilon=kh\ll 1$ since $k\in{\cal B}$ and $h\ll1$. We define the following integral operators:
\begin{align}
  \label{eq:int:S} [{\cal S}\phi](x) =& \int_{\Gamma_h} \frac{e^{\bi k |x-y|}}{2\pi|x-y|}\phi(y)dS(y),\quad x\in\Gamma_h,\\
  \label{eq:int:S0}
  [{\cal S}_0\phi](x) =& \int_{\Gamma_1}\frac{1}{2\pi|x-y|}\phi(y)dS(y),\quad x\in\Gamma_1,\\
  \label{eq:int:R}
  [{\cal R}_0\phi](x) =& \int_{\Gamma_1}\frac{e^{\bi \epsilon|x-y|}-1-\bi\epsilon|x-y|}{2\pi\epsilon^2|x-y|}\phi(y)dS(y),\quad x\in\Gamma_1,
\end{align}
where $\Gamma_1=G_1\times\{0\}$ is the
aperture surface $\Gamma_h$ when $h=1$. Some important properties of the above
integral operators are listed below.
\begin{mylemma}
  \label{lem:propSR}
  For any $k\in{\cal B}$ and $h\ll1$, we have:
  \begin{itemize}
    \item[1.] ${\cal S}$ can be uniquely extended as a bounded operator from
      $\tilde{H}^{-1/2}(\Gamma_h)$ to $H^{1/2}(\Gamma_h)$;
    \item[2.] ${\cal S}_0$ can be uniquely extended as a bounded operator from
      $\tilde{H}^{-1/2}(\Gamma_1)$ to $H^{1/2}(\Gamma_1)$, and ${\cal S}_0$ is
      positive and bounded below on $\tilde{H}^{-1/2}(\Gamma_1)$, i.e., for any
      $\phi\in \tilde{H}^{-1/2}(\Gamma_1)$,
      \[
        <{\cal S}_0\phi,\phi>_{1/2,-1/2}\geq C||\phi||_{\tilde{H}^{-1/2}(\Gamma_1)},
      \]
      for some positive constant $C>0$.
    \item[3.] ${\cal R}_0$ can be uniquely extended as a uniformly bounded
      operator from $\tilde{H}^{-1/2}(\Gamma_1)$ to $H^{1/2}(\Gamma_1)$, i.e.,
      \[
        ||{\cal R}_0||\leq C,
      \]
      for some constant $C>0$ independent of $\epsilon$.
  \end{itemize}
  \begin{proof}
    Choosing any bounded and closed Lipschitz surface
    $\Gamma_c\subset\overline{\mathbb{R}_+^3}$ that contains $\Gamma_1$, we can
    extend ${\cal S}_0$ as for any $\phi\in \tilde{H}^{-1/2}(\Gamma_1)\subset
    H^{-1/2}(\Gamma_c)$,
    \[
      {\cal S}_0\phi = (2{\cal S}_c\phi)|_{\Gamma_1},
    \]
    where
    \[
      {\cal S}_c\phi = \int_{\Gamma_c}G_0(x,y)\phi(y) dS(y),
    \]
    with the three dimensional fundamental solution
    $G_0(x,y)=\frac{1}{4\pi|x-y|}$ of Laplacian as the kernel.
    According to \cite[Thm. 7.6 \& Cor. 8.13]{mcl00}, it is clear that ${\cal S}_c$ is
    bounded from $H^{-1/2}(\Gamma_c)$ to $H^{1/2}(\Gamma_c)$, and satisfies for
    any $\psi\in H^{-1/2}(\Gamma_c)$,
    \[
      <{\cal S}_c\psi,\psi>_{H^{1/2}(\Gamma_c),H^{-1/2}(\Gamma_c)}\geq C||\psi||_{H^{-1/2}(\Gamma_c)}/2,
    \]
    for some constant $C>0$. Thus, for any $\phi\in \tilde{H}^{-1/2}(\Gamma_1)$,
    ${\cal S}_0\phi\in H^{1/2}(\Gamma_1)$ so that
    \[
      <{\cal S}_0\phi,\phi>_{H^{1/2}(\Gamma_1),H^{-1/2}(\Gamma_1)} = <2{\cal S}_c\phi,\phi>_{H^{1/2}(\Gamma_c),H^{-1/2}(\Gamma_c)}\geq C||\phi||_{\tilde{H}^{-1/2}(\Gamma_1)}.
    \]
    The mapping property of ${\cal S}$ on $\Gamma_h$ can be similarly obtained.
    Now we prove the uniform boundedness of ${\cal R}_0$. Since the kernel
    function of ${\cal R}_0$ in (\ref{eq:int:R}) and its gradient with respect
    to $x$ are uniformly bounded for all $k\in{\cal B}$ and $h\ll 1$, we easily
    conclude that ${\cal R}_0$ is uniformly bounded from $L^2(\Gamma_1)$ to
    $H^{1}(\Gamma_1)$. The self-adjointness of ${\cal R}_0$ and interpolation
    theory then indicate the uniform boundedness of ${\cal R}_0$ mapping from
    $\tilde{H}^{-1/2}(\Gamma_1)$ to $H^{1/2}(\Gamma_1)$.
   \end{proof}
\end{mylemma}
In the upper half plane $\mathbb{R}_+^3$, the Neumann data
$\partial_{\nu}u(x)=-\partial_{x_3}u(x)$ in $\tilde{H}^{-1/2}(\Gamma_h)$
uniquely determines the following outgoing solution
\[
  u(x) = -\int_{\Gamma_h} \frac{e^{\bi k |x-y|}}{2\pi|x-y|}\partial_{x_3} u(y)dS(y)\in H^{1}_{\rm loc}(\mathbb{R}_+^3),\quad{\rm on}\quad\mathbb{R}_+^3,
\]
with trace $u|_{\Gamma_h} = -{\cal S}\partial_{x_3}u\in H^{1/2}(\Gamma_h)$ by Lemma~\ref{lem:propSR}.
By the continuity of $\partial_{\nu}u$ on $\Gamma_h$, to ensure that $u\in H^1_{\rm loc}(\Omega_h^+)$, we require that
\[
  u|_{\Gamma_h} = -{\cal S}\partial_{x_3} u =
  \sum_{m=0}^{\infty}b_{m}\phi_{m}(x;h)[e^{\bi s_{m} l} + 1]\in H^{1/2}(\Gamma_h),
\]
which is equivalent to for any nonnegative integers $m$,
\[
  (-{\cal S}\partial_{x_3} u,\phi_{m}(x_1,x_2;h))_{L^2(\Gamma_h)} =
  b_{m}[e^{\bi s_{m}l}+1].
\]
By (\ref{eq:normder}), the above equations can be rewritten as the following
infinite number of linear equations in terms of unknowns
$\{b_{0},\{a_{m}\}\}\in\ell^2$, i.e.,
\begin{align}
  \label{eq:b0}
  b_{0}[e^{\bi s_{0}l}+1] =&- b_{0}\bi s_{0}h[e^{\bi s_{0}l} - 1]d_{00} - \sum_{m'=1}^{\infty}a_{m'}\lambda_{m'}^{-1/4}\bi s_{m'}h[e^{\bi s_{m'}l} - 1]d_{m'0},\\
  \label{eq:am}
  \lambda_m^{-1/4}a_{m}[e^{\bi s_{m}l}+1] =&- b_{0}\bi s_{0}h[e^{\bi s_{0}l} - 1]d_{0m}
- \sum_{m'=1}^{\infty}a_{m'}\lambda_{m'}^{-1/4}\bi s_{m'}h[e^{\bi s_{m'}l} - 1]d_{m'm},
\end{align}
where we have set for integers $m',m\geq 0$ that
\begin{equation}
  \label{eq:dm'm}
d_{m'm}=h^{-1}({\cal S}\phi_{m'}(\cdot;h),\phi_{m}(\cdot;h))_{L^2(\Gamma_h)}.
\end{equation}
We have the following properties regarding asymptotics of $d_{m'm}$ as
$\epsilon\ll 1$.
\begin{mylemma}
  \label{lem:dm'm}
  For $\epsilon\ll1$ and for nonnegative integers $m,m'$, $d_{m'm}$
  asymptotically behaves as following:
\begin{align}
  \label{eq:dm'm}
  d_{m'm} = \left\{
  \begin{array}{ll}
    ({\cal S}_01,1)_{L^2(\Gamma_1)}+ \frac{\bi \epsilon}{2\pi} + {\cal O}(\epsilon^2), & m=m'=0;\\
    ({\cal S}_{0}\phi_{m'},\phi_{m})_{L^2(\Gamma_1)} + \epsilon^2({\cal R}_0\phi_{m'},\phi_{m})_{L^2(\Gamma_1)}, & {\rm otherwise},
  \end{array}
  \right.
\end{align}
where $\phi_{m}$ is short for $\phi_{m}(\cdot;h)$ when $h=1$.
  \begin{proof}
    By rescaling,
    \begin{align*}
      d_{m'm} =& \frac{1}{2\pi}\int_{\Gamma_1}\int_{\Gamma_1}\frac{e^{\bi \epsilon|x-y|}}{|x-y|}\phi_{m'}(y;1)dS(y)\bar{\phi}_{m}(x;1)dS(x),\\
      =&({\cal S}_0\phi_{m'},\phi_{m})_{L^2(\Gamma_1)} +\frac{\bi\epsilon }{2\pi}\int_{\Gamma_1}\int_{\Gamma_1}\phi_{m'}(y;1)dy_1dy_2\bar{\phi}_{m}(x;1)dx_1dx_2 +\epsilon^2({\cal R}_0\phi_{m'},\phi_{m})_{L^2(\Gamma_1)}.
    \end{align*}
    It is clear that the second term of the r.h.s is nonzero only when
    $m=m'=0$.
  \end{proof}
\end{mylemma}

Now set for $m>0$ and $m'>0$ that
\begin{align}
  c_{00} =& -\bi \epsilon d_{00},\\
  c_{m'0} =& -\lambda_{m'}^{-1/4}\bi s_{m'}h[e^{\bi s_{m'}l} - 1] d_{m'0},\\
  c_{0m} =&- \lambda_m^{1/4}\frac{\bi \epsilon}{e^{\bi s_{m}l}+1}d_{0m},\\
  c_{m'm} =&-\lambda_{m'}^{-1/4}\lambda_m^{1/4}\bi s_{m'}h\frac{e^{\bi s_{m'}l} - 1}{e^{\bi s_{m}l}+1}d_{m'm},
\end{align}
we could rewrite the previous linear equations (\ref{eq:b0}) and (\ref{eq:am})
more compactly in the following form,
\begin{align}
  \label{eq:cpb}
  b_{0}[e^{\bi s_{0}l}+1] =&b_{0}[e^{\bi s_{0}l} - 1]c_{00} + <\{a_{m'}\},\{c_{m'0}\}>_{\ell^2},\\
  \label{eq:cpam}
  \{a_{m}\}=&b_{0}[e^{\bi s_{0}l} - 1]\{c_{0m}\}+ {\cal A}_h\{a_{m}\},
\end{align}
where the operator ${\cal A}_h$ is defined as: for any $\{f_{m}\}_{m=1}^{\infty}\in\ell^2$,
\begin{align}
  \label{eq:Ah}
  {\cal A}_h\{f_{m}\}_{m=1}^{\infty} = \{\sum_{m=1}^{\infty} c_{m'm}f_{m'}\}_{m=1}^{\infty},
\end{align}
According to Lemma~\ref{lem:dm'm}, we have the following properties.
\begin{mylemma}
  \label{lem:cm'mA}
  For $h\ll 1$ and $k\in{\cal B}$:
  \begin{itemize}
  \item[1.]
      \[
        c_{00} = -({\cal S}_01,1)_{L^2(\Gamma_1)}\epsilon\bi+ \frac{\epsilon^2}{2\pi} + {\cal O}(\epsilon^3);
      \]
    \item[2.] When $m>0$,
      \[
        c_{0m}= -\bi\epsilon({\cal S}_{0}1,\lambda_m^{1/4}\phi_{m})_{L^2(\Gamma_1)} -\bi\epsilon^3({\cal R}_01,\lambda_m^{1/4}\phi_{m})_{L^2(\Gamma_1)},
      \]
      and $\{c_{0m}\}_{m=1}^{\infty}\in \ell^2$;
    \item[3.] When $m'>0$,
      \begin{align*}
        c_{m'0}=& -({\cal S}_{0} 1,\lambda_{m'}^{1/4}\phi_{m'})_{L^2(\Gamma_1)}- \epsilon^2({\cal R}_0 1, \lambda_{m'}^{1/4}\phi_{m'n'})_{L^2(\Gamma_1)} + ({\cal S}_0 1, \lambda_{m'}^{-3/4}\phi_{m'n'}){\cal O}(\epsilon^2),
      \end{align*}
      and $\{c_{m'0}\}_{m'=1}^{\infty}\in \ell^2$;
    \item[4.] When $m',m>0$,
      \begin{align*}
        c_{m'm} =& -({\cal S}_0\lambda_{m'}^{1/4}\phi_{m'},\lambda_m^{1/4}\phi_{m})_{L^2(\Gamma_1)}-\epsilon^2({\cal R}_0\lambda_{m'}^{1/4}\phi_{m'},\lambda_m\phi_{m})_{L^2(\Gamma_1)} +({\cal S}_0\lambda_{m'}^{-3/4}\phi_{m'},\lambda_m^{1/4}\phi_{m})_{L^2(\Gamma_1)}{\cal O}(\epsilon^2),
      \end{align*}
      and the operator ${\cal A}_h$ defined by
      $\{c_{m'm}\}_{m,m'>0}$ is bounded from $\ell^2$ to
      $\ell^2$ and can be decomposed as
      \begin{align}
        {\cal A}_h = {\cal P} + \epsilon^2{\cal Q}_h,
      \end{align}
      where ${\cal P}$ is defined as
      \[
        {\cal P}\{f_{m}\}_{m>0} = \left\{-\sum_{m'>0}({\cal S}_0\lambda_{m'}^{1/4}\phi_{m'},\lambda_m^{1/4}\phi_{m})_{L^2(\Gamma_1)}f_{m}\right\}_{m>0},
      \]
      and ${\cal Q}_h = \epsilon^{-2}({\cal A}_h - {\cal P})$. Both ${\cal P}$
      and ${\cal Q}_h$ are uniformly bounded from $\ell^2$ to $\ell^2$ for all
      $k\in{\cal B}$ and $h\ll 1$.
  \end{itemize}
  \begin{proof}
    The asymptotic behaviors of $\{c_{m'm}\}_{m,m'>0}$ are
    trivial by Lemma~\ref{lem:dm'm}. As for the other properties, we here
    only show that the leading terms of $\{c_{m'm}\}$ satisfy those
    properties as the high-order terms can be analyzed similarly. For any
    $\{f_{m}\}_{m>0}\in\ell^2$ and $\{g_{m}\}_{m>0}\in\ell^2$, the
    following two functions
   \begin{align*}
     f = \sum_{m>0}f_{m}\lambda_m^{1/4}\phi_{m},\ {\rm and}\ g = \sum_{m>0}f_{m}\lambda_m^{1/4}\phi_{m},
   \end{align*}
   are in $H^{-1/2}(\Gamma_1)$. We thus have
   \begin{align*}
     &|\sum_{m>0}({\cal S}_{0}1,\lambda_m^{1/4}\phi_{m})_{L^2(\Gamma_1)}\bar{f}_{m}| \\
     =&|<-{\cal S}_0 1, f>_{1/2,-1/2}| \leq ||{\cal S}_0 1||_{H^{1/2}(\Gamma_1)}||f||_{H^{-1/2}(\Gamma_1)} \leq ||{\cal S}_0 1||_{H^{1/2}(\Gamma_1)}||\{f_{m}\}||_{\ell^2},
   \end{align*}
  indicating that $\{c_{0m}\}\in \ell^2$ and
  $\{c_{m'0}\}\in \ell^2$. Moreover,
  \begin{align*}
    &\left|\sum_{m>0}\sum_{m'>0}-({\cal S}_0\lambda_{m'}^{1/4}\phi_{m'},\lambda_m^{1/4}\phi_{m})_{L^2(\Gamma_1)}f_{m'}g_{m}\right| \\
    =& |<{\cal S}_0f,g>_{1/2,-1/2}| \leq ||{\cal S}_0||\cdot ||\{f_{m}\}||_{\ell^2}||g_{m}||_{\ell^2},
  \end{align*}
  which implies the mapping property of ${\cal A}_h$ as well as the
  boundedness.
\end{proof}
\end{mylemma}
\subsubsection{Invertibility of ${\cal I}-{\cal A}_h$}
To prove ${\cal I}-{\cal A}_h$ has a bounded inverse for $h\ll 1$, we no longer
justify the diagonal dominance of the infinite-dimensional matrix
$\{\delta_{m'm}-c_{m'm}\}_{m',m=0}^{\infty}$ as was done in \cite{zhlu21}, since now
$c_{m'm}$ is challenging to accurately approximate, although we conjecture that
this property remains true\footnote[1]{Such a property has been numerically
  verified for a rectangular hole.}. To resolve this issue, we convert the study
of the invertibility of ${\cal I}-{\cal A}_h$ to the study of the well-posedness
of a closely related boundary-value problem in a semi-infinite cylinder
$\Omega^-=G_1\times \mathbb{R}^-$, as was done in \cite{bontri10}. Nevertheless,
compared with the two-dimensional proof in \cite{bontri10}, our three-dimensional
proof is simpler, more versatile, and more straightforward due to the following
aspects: (1) we make no use of the Green function of $\Omega^-$, which is
complicated; (2) it is not necessary to introduce a Dirichlet-to-Neumann map to
truncate the unbounded cylinder $\Omega^-$; (3) the proof does not use any
particular property regarding the shape of $\Omega^-$ (e.g. a circular hole in
\cite{liazou20} or a rectangular hole).

 As seen in
Lemma~\ref{lem:cm'mA}, ${\cal P}$ can be regarded as the limit of ${\cal A}_h$
as $h\to 0$ so that by Neumann series, we could see that ${\cal I}-{\cal A}_h$
is invertible if ${\cal I}-{\cal P}$ is invertible. The open mapping theorem
implies the following lemma.
\begin{mylemma}
  \label{lem:P1}
  For $k\in{\cal B}$ and $h\ll1$, the operator ${\cal I} - {\cal P}:\ell^2\to
  \ell^2$ has a bounded inverse, if for any
  $\{\alpha_{m}\}_{m>0}\in\ell^2$, the following problem
\begin{align*}
  {\rm (P1):}\quad \{a_{m}\}_{m>0}= {\cal P}\{a_{m}\}_{m>0} + \{\alpha_{m}\}_{m>0},
\end{align*}
has a unique solution $\{a_{m}\}_{m>0}\in\ell^2$.
\end{mylemma}
For $\{\alpha_{m}\}_{m>0}\in\ell^2$, let
\begin{align}
  \label{eq:beta}
  \beta_{m}=&\alpha_{m}\lambda_{m}^{-1/4},\quad m>0,\\
  \label{eq:f}
  f(x_1,x_2) =& \sum_{m>0}\beta_{m}\phi_{m}(x),\quad |x_i|\leq 1/2,i=1,2.
\end{align}
Let $H_0^{1/2}(\Gamma_1)=\{\phi\in H^{1/2}(\Gamma_1):(\phi,1)_{L^2(\Gamma_1)}=0\}$ and
$\tilde{H}_0^{-1/2}(\Gamma_1)=\{\phi\in \tilde{H}^{-1/2}(\Gamma_1):<\phi,1>_{-1/2,1/2}=0\}$. It can be seen that
$f\in H_0^{1/2}(\Gamma_1)$. Now, consider the following problem:
\begin{align*}
  {\rm (P2):}\quad\left\{
  \begin{array}{ll}
    \Delta u = 0,&{\rm on}\quad \Omega^-=G_1\times\mathbb{R}^-,\\
    \partial_{\nu} u = 0,&{\rm on}\quad \partial\Omega^-\backslash\overline{\Gamma_1},\\
    u +{\cal S}_0\partial_{\nu}u =f +<S_0\partial_{\nu}u-f,1>_{1/2,-1/2},&{\rm on}\quad \Gamma_1,\\
    %<\partial_{\nu} u, 1>_{\Gamma_A} = \int_{\Gamma_A}\partial_\nu u\cdot 1 dS = 0.
  \end{array}
  \right.
\end{align*}
We have the following lemma.
\begin{mylemma}
  \label{lem:P1P2}
  For any $\{\alpha_{m}\}_{m>0}\in\ell^2$, problem (P1) has a unique
  solution $\{a_{m}\}_{m>0}\in\ell^2$ iff problem (P2) has a unique solution
  $u\in H^{1}(\Omega^-)$ for any $f\in H_0^{1/2}(\Gamma_1)$ defined in (\ref{eq:f}).
  \begin{proof}
    Given $\{\alpha_{m}\}_{m>0}\in\ell^2$, suppose (P1) has a
    solution $\{a_{m}\}_{m>0}\in\ell^2$. Setting
    $b_{m}=a_{m}\lambda_m^{-1/4}$ for $m>0$, we now claim that
    \[
      u(x) = \sum_{m=1}^{\infty}b_{m}\phi_{m}(x)e^{\lambda_m^{1/2} x_3}\in H^{1}(\Omega^-),
    \]
    and solves (P2).
    Clearly,
    \begin{align*}
      ||u||_{H^1(\Omega^-)} \eqsim & \int_{-\infty}^0dx_3\int_{G_1}|u(x_1,x_2,x_3)|^2 + |\partial_{x_1}u(x_1,x_2,x_3)|^2 + |\partial_{x_2}u(x_1,x_2,x_3)|^2 dx_1dx_2\\
      +&\int_{-\infty}^0dx_3\int_{G_1}|\partial_{x_3}u(x_1,x_2,x_3)|^2dx_1dx_2\\
      =& \int_{-\infty}^0dx_3 ||u(\cdot,x_3)||^2_{H^1(G_1)}dx_1dx_2 + \int_{-\infty}^0dx_3||\partial_{x_3}(\cdot,x_3)||^2_{L^2(G_1)}\\
      \eqsim &\int_{-\infty}^0 \sum_{m>0}^{\infty}(1+\lambda_m)|b_{m}|^2e^{2\lambda_m^{1/2} x_3}dx_3 +  \int_{-\infty}^0\sum_{m>0}\lambda_m|b_m|^2e^{2\lambda_m^{1/2}x_3}dx_3\\
      =&\sum_{m>0}|a_m|^2\frac{1+2\lambda_m}{2\lambda_m}<\infty.
    \end{align*}
    Now, we verify that $u$ solves (P2). It is clear that $\Delta u = 0$ on
    $\Omega^-$ in the distributional sense so that on $\Gamma_1$,
    \[
      \partial_{\nu}u = \sum_{m>0} b_{m}\phi_{m}(x)\lambda_{m}^{1/2}\in ( H^{1/2}(\Gamma_1))'=\tilde{H}^{-1/2}(\Gamma_1),
    \]
    since $\partial_{\nu}u=0$ on
    $\partial\Omega^-\backslash\overline{\Gamma_1}$. Clearly, $<\partial_{\nu}u,
    1>_{-1/2,1/2}=0$ so that $\partial_{\nu}u\in \tilde{H}^{-1/2}_0(\Gamma_1)$. Thus, it
    suffices to prove that
    \[
      u +{\cal S}_0\partial_{\nu}u =f + <S_0\partial_{\nu}u-f,1>_{1/2,-1/2},
    \]
    in $H^{1/2}(\Gamma_1)$. In fact, for any positive integers $m'>0$,
    \begin{align*}
      <u +{\cal S}_0\partial_{\nu}u - f,\phi_{m'}>_{1/2,-1/2} =& b_{m'}  +\sum_{m>0}b_{m}\lambda_{m}^{1/2}<{\cal S}_0\phi_{m},\phi_{m'}>_{1/2,-1/2}-\beta_{m'} =0.
    \end{align*}
    Conversely, suppose $u\in H^1(\Omega^-)$ solves (P2). We choose for $m>0$,
    \begin{align*}
      b_{m} = (u|_{\Gamma_1}, \phi_m)_{L^2(\Gamma_1)},
    \end{align*}
    so that $\{a_{mn}=b_{m}\lambda^{-1/4}\}\in\ell^2$ since $u|_{\Gamma_1}\in
    H^{1/2}(\Gamma_1)$. Now, take
    \[
      u^* = \sum_{m>0}b_{m}\phi_{m}(x)e^{\lambda_m^{1/2}x_3}\in H^{1}(\Omega^-),
    \]
     and we claim that $u-u^*=0$ on $\Omega^-$. In fact, $\phi = u-u^*\in H^1(\Omega^-)$ solves
     \begin{equation*}
       \left\{
         \begin{array}{ll}
           \Delta \phi = 0,&{\rm on}\quad\Omega^-,\\
           \phi = 0,&{\rm on}\quad\Gamma_1,\\
           \partial_{\nu}\phi = 0,&{\rm on}\quad\partial\Omega^-\backslash\overline{\Gamma_1}.
         \end{array}
       \right.
     \end{equation*}
       Thus, testing the governing equation with $\phi$ itself yields
       \[
         -||\nabla\phi||^2_{\Omega^-}=0,
       \]
       so that $\phi$ is constant on $\Omega^-$. But $\phi|_{\Gamma_1}=0$
       implies that $\phi=0$. Consequently, following a similar argument as
       before, we could verify that $\{a_{m}\}_{m>0}$ solves (P1).
  \end{proof}
\end{mylemma}
To prove that (P2) has a unique solution, we make use of the method of
variational formulation. Let
\[
  V = H^1(\Omega^-)\times \tilde{H}^{-1/2}(\Gamma_1),
\]
be equipped with the natural cross-product norm, and let
$a: V\times V\to \mathbb{C}$ be defined as:
\begin{align}
  a((u,\phi),(v,\psi)) =& (\nabla u,\nabla v)_{L^2(\Omega^-)} - <\phi, v>_{-1/2,1/2}+ <u, \psi>_{1/2,-1/2}\\
  &+ <{\cal S}_0\phi,\psi-<\psi,1>_{-1/2,1/2}>_{1/2,-1/2},
\end{align}
for any $(u,\phi),(v,\psi)\in V$. Such a formulation of $a$ can be obtained by
testing the first equation of (P2) with $v$ and the third equation with $\psi$.
Then, (P2) is equivalent to the
following variational problem: Find $(u,\phi)\in V$, s.t.,
\begin{align}
  {\rm (P3)}:\quad a((u,\phi),(v,\psi)) = <f,\psi-<\psi,1>_{-1/2,1/2}>_{1/2,-1/2},
\end{align}
for all $(v,\psi)\in V$. Though Lemma~\ref{lem:P1P2} requires that $f\in
H^{1/2}_0(\Gamma_1)$, it turns out that $f\in H^{1/2}(\Gamma_1)$ is also allowed
as illustrated in the following theorem.
\begin{mytheorem}
  \label{thm:P3}
  For any $f\in H^{1/2}(\Gamma_1)$, the variational problem (P3) has a unique solution.
  \begin{proof}
    We first prove that $\phi\in \tilde{H}_{0}^{-1/2}(\Gamma_1)$, i.e.,
    $<\phi,1>_{-1/2,1/2}=0$. Let  \[
    v_n(x) = \left\{
      \begin{array}{ll}
      1,& \{x\in \Omega^-:x_3\in(-n,0)\};\\
      e^{-(x_3+n)^2},& {\rm otherwise},
      \end{array}
    \right.
  \]
  and $\psi=0$ so that
  \[
    \int_{\{x\in \Omega^-:x_3\leq -n\}}\nabla u \nabla \overline{v_n} dx =
    <\phi,1>_{-1/2,1/2}.
  \]
  Letting $n\to \infty$ yields that $<\phi,1>_{-1/2,1/2}=0$. Now, by Lemma~\ref{lem:propSR},
  \[
    {\rm Re}[a((u,\phi),(u,\phi))] = ||\nabla u||^2_{L^2(\Omega^-)} + <{\cal
      S}_0\phi,\phi>_{1/2,-1/2}\geq ||u||^2_{H^1(\Omega^-)}+ C||\phi||^2_{H^{-1/2}(\Gamma_1)}
    - ||u||^2_{L^2(\Omega^-)},
  \]
implies that the bilinear functional $a$ defines a Fredholm operator of index
zero so that we only need to show the uniqueness. Now suppose $f=0$, then the above equation in fact implies
\[
  0 = ||\nabla u||^2_{L^2(\Omega^-)} + <{\cal
      S}_0\phi,\phi>_{1/2,-1/2}\geq ||\nabla u||^2_{L^2(\Omega^-)} + C||\phi||^2_{H^{-1/2}(\Gamma_1)},
\]
so that $u$ must be a constant in $\Omega^-$ and $\phi\equiv 0$ in
$H^{-1/2}(\Gamma_1)$. Choosing $v = 0$ and $\psi = 1$, we get
\[
 <u,1>_{1/2,-1/2} = 0,
\]
which implies that $u\equiv 0$ in $\Omega^-$. Finally, the proof is concluded from
that the r.h.s of (P3) defines a bounded functional in $V^*$ for any $f\in
H^{1/2}(\Gamma_1)$.
  \end{proof}
\end{mytheorem}
Combining Theorem~\ref{thm:P3}, Lemmas~\ref{lem:P1} and \ref{lem:P1P2} yields
the desired result.
\begin{mytheorem}
  \label{thm:PAbi}
  For any $k\in{\cal B}$ and $h\ll1$, both ${\cal I}-{\cal P}$ and ${\cal
    I}-{\cal A}_h$ have bounded inverses. In fact,
  \[
    ||({\cal I}-{\cal A}_h)^{-1} - ({\cal I}-{\cal P})^{-1}|| = {\cal
      O}(\epsilon^2),\quad \epsilon\ll1.
  \]
  \begin{proof}
   It is clear that the two operators have bounded inverses. Now, we prove the
   estimate. By Lemma~\ref{lem:cm'mA},
   \begin{align*}
     ({\cal I}-{\cal A}_h)^{-1} = ({\cal I}-{\cal P} + \epsilon^2{\cal Q}_h)^{-1}
     = ({\cal I}-{\cal P})^{-1}({\cal I}+\epsilon^2{\cal Q}_h({\cal I}-{\cal P})^{-1})^{-1},
   \end{align*}
   so that based on Neumann series,
   \begin{align*}
     ||({\cal I}-{\cal A}_h)^{-1} - ({\cal I}-{\cal P})^{-1}|| =&||({\cal I}-{\cal P})^{-1}(\sum_{n=1}^{\infty}(-\epsilon^2{\cal Q}_h({\cal I}-{\cal P})^{-1})^{n}||\\
     \leq& \epsilon^2||({\cal I}-{\cal P})^{-1}||^2\cdot ||{\cal Q}_h||\cdot\frac{1}{1-\epsilon^2||{\cal Q}_h({\cal I}-{\cal P})^{-1}||}={\cal O}(\epsilon^2).
   \end{align*}
  \end{proof}
\end{mytheorem}
  Finally, since ${\cal S}_0$ is strictly positive definite, ${\cal I}-{\cal P}$
  as well as its inverse is also strictly positive definite in the sense that:
  for any $\{f_m\}\in \ell^2$,
  \begin{equation}
    \label{eq:spd:P}
    <({\cal I}-{\cal P})^{-1}f_m,f_m>_{\ell^2} \geq C||\{f_m\}||_{\ell^2},
  \end{equation}
  for some positive constant $C$.

\subsubsection{Resonance Frequencies}
Based on Theorem~\ref{thm:PAbi}, (\ref{eq:cpb}) and (\ref{eq:cpam}) are reduced
to the following single equation for the unknown $b_{00}$.
\begin{align}
  \left[  (e^{\bi k l}+1)-(e^{\bi k l} - 1)\left(c_{00} + <({\cal I}-{\cal A}_h)^{-1}\{c_{0m}\},\{c_{m0}\}>_{\ell^2}  \right)\right]b_{00} = 0.
\end{align}
Now, based on Lemma~\ref{lem:cm'mA} and Theorem~\ref{thm:PAbi}, we obtain our
first main result.
\begin{mytheorem}
  \label{thm:evenres}
  For any width $h\ll 1$, the governing equations (\ref{eq:gov:e}) and
  (\ref{eq:bc:e}) possess nonzero solutions for $k\in{\cal B}$, if and only if
  the following nonlinear equation of $k$
  \begin{equation}
    \label{eq:gov:ke}
    (e^{\bi k l}+1)-(e^{\bi k l} - 1)\left(c_{00} + <({\cal I}-{\cal A}_h)^{-1}\{c_{0m}\},\{c_{m0}\}>_{\ell^2} \right)= 0,
  \end{equation}
  has solutions in ${\cal B}$. In fact, these solutions (the so-called resonance
  frequencies) are
  \begin{equation}
    \label{eq:ke:asy}
    kl = k_{m,e} -2\bi \Pi(\epsilon_{m,e}) + 2k_{m,e}^{-1}\Pi(\epsilon_{m,e})(\pi^{-1}\epsilon_{m,e} + 2\Pi(\epsilon_{m,e}) ) + {\cal O}(\epsilon_{m,e}^3), m=1,2,\cdots.
  \end{equation}
  where $k_{m,e} = (2m-1)\pi$ is a Fabry-P\'erot frequency,
  $\epsilon_{m,e}=k_{m,e}h\ll 1$, and
  \begin{align}
    \label{eq:alpha}
    \alpha =& <({\cal I}-{\cal P})^{-1}\{({\cal S}_01,\lambda_m^{1/4}\phi_{m})_{L^2(\Gamma_1)}\},\{({\cal S}_01,\lambda_m^{1/4}\phi_{m})_{L^2(\Gamma_1)}\}>_{\ell^2}\geq 0,\\
    \label{eq:Pi}
    \Pi(\epsilon) =& \left[  -({\cal S}_01,1)_{L^2(\Gamma_1)}+\pi\alpha\right]\epsilon\bi+ \frac{\epsilon^2}{2\pi}.
  \end{align}
  \begin{proof}
    For $h\ll 1$, $\epsilon\ll1$ so that by Theorem~\ref{thm:PAbi} and
    Lemma~\ref{lem:cm'mA}, equation (\ref{eq:gov:ke}) can be reduced to
    \[
      e^{\bi kl}+1 = (e^{\bi kl}-1)\Pi(\epsilon) + {\cal O}(\epsilon^3),
    \]
    which is equivalent to
    \[
      e^{\bi kl} + 1 = -\frac{2\Pi(\epsilon)}{1-\Pi(\epsilon)} + {\cal O}(\epsilon^3).
    \]
    Here, by (\ref{eq:spd:P}), it can be easily shown that $\alpha>0$.

    As the right-hand side approaches $0$ as $\epsilon\to 0$, we see that the
 resonance frequencies must satisfy: for some $m=1,\cdots$, $\delta_{m,e}:= k l - k_{m,e} = o(1)$. Thus, $\epsilon - \epsilon_{m,e} =h\delta_{m,e}$, as $h \to 0^+$. Therefore, we have
 \[
   e^{\bi \delta_{m,e}} - 1 = \frac{2\Pi(\epsilon)}{1-\Pi(\epsilon)} + {\cal O}(\epsilon^3),
 \]
 so that by Taylor's expansion of $\log(1+2x/(1-x))$ at $x=0$,
 \begin{align*}
   \delta_{m,e} =& -\bi\log\left[  1 +\frac{2\Pi(\epsilon)}{1-\Pi(\epsilon)} + {\cal O}(\epsilon^3)\right]=-2\bi\Pi(\epsilon) + {\cal O}(\epsilon^3).
   %=&-\bi\left[ \frac{\Pi(\epsilon)}{1-\Pi(\epsilon)/2} -\frac{\Pi^2(\epsilon)}{2(1-\Pi(\epsilon)/2)^2} + \frac{\Pi^3(\epsilon)}{3(1-\Pi(\epsilon)/2)^3} \right] + {\cal O}(\epsilon^3\log\epsilon)\\
   %=&-\bi\left[\Pi(\epsilon) + \frac{1}{12}\Pi^3(\epsilon) \right] + {\cal O}(\epsilon^3\log\epsilon).
 \end{align*}
 Thus,
 \[
   \delta_{m,e}\eqsim 2\left[  -\frac{2}{\pi}\log(\sqrt{2}+1)+\frac{2}{3\pi}(\sqrt{2}-1)+\pi\alpha\right]\epsilon_{m,e},
 \]
 Based on the definition of $\Pi$, we get
 \begin{align*}
   \Pi(\epsilon)-\Pi(\epsilon_{m,e}) =& k_{m,e}^{-1}\delta_{m,e}\left[  -\frac{2}{\pi}\log(\sqrt{2}+1)+\frac{2}{3\pi}(\sqrt{2}-1)+\pi\alpha\right]\epsilon_{m,e}\bi \\
   &+ \frac{k_{m,e}^{-1}\epsilon_{m,e}}{2\pi}\delta_{m,e}(k_{m,e}^{-1}\epsilon_{m,e}\delta_{m,e}+2\epsilon_{m,e})+ {\cal O}(\epsilon_{m,e}^3)\\
   =&k_{m,e}^{-1}\delta_{m,e}\left(  \Pi(\epsilon_{m,e}) + \frac{\epsilon_{m,e}}{2\pi} \right) +\frac{(k_{m,e}^{-1}\epsilon_{m,e})^2}{2\pi}\delta_{m,e}^2 + {\cal O}(\epsilon_{m,e}^3).
 \end{align*}
 Thus
 \begin{align*}
   \delta_{m,e} =&-2\bi \Pi(\epsilon_{m,e}) -2\bi k_{m,e}^{-1}\delta_{m,e}\left(  \Pi(\epsilon_{m,e})+\frac{\epsilon_{m,e}}{2\pi} \right)+{\cal O}(\epsilon_{m,e}^3),
 \end{align*}
 so that
\begin{align*}
  \delta_{m,e} %=& -\frac{2C}{B}\frac{1}{1+\sqrt{1-\frac{4AC}{B^2}}}
  %=&-\frac{2C}{B}\left[ \frac{1}{2} + \frac{AC}{2B^2} \right] + {\cal O}(\epsilon_{m,e}^5\log^3\epsilon_{m,e})\\
  %=-\frac{C}{B} - \frac{AC^2}{B^3} + {\cal O}(\epsilon_{m,e}^5)\\
  =&\frac{2\bi \Pi(\epsilon_{m,e})}{-\frac{\bi}{\pi} k_{m,e}^{-1}\epsilon_{m,e} - 2\bi k_{m,e}^{-1}\Pi(\epsilon_{m,e}) - 1} + {\cal O}(\epsilon_{m,e}^3)\\
  =&-2\bi \Pi(\epsilon_{m,e}) + 2\Pi(\epsilon_{m,e})(\pi^{-1}k_{m,e}^{-1}\epsilon_{m,e} + 2k_{m,e}^{-1}\Pi(\epsilon_{m,e}) ) + {\cal O}(\epsilon_{m,e}^3).
  %=&-\bi \Pi(\epsilon_{m,e})\left[1 + \left(  \frac{2}{\pi} k_{m,e}^{-1}\epsilon_{m,e} - \bi k_{m,e}^{-1}\Pi(\epsilon_{m,e})\right) + \left(  \frac{2}{\pi} k_{m,e}^{-1}\epsilon_{m,e} - \bi k_{m,e}^{-1}\Pi(\epsilon_{m,e})\right)^2  \right] \\
  %&-\frac{\bi}{12}\Pi^3(\epsilon_{m,e})-\pi^{-1}k_{m,e}^{-2}\epsilon_{m,e}\Pi^2(\epsilon_{m,e}) + {\cal O}(\epsilon_{m,e}^3\log\epsilon_{m,e})\\
  %=&-\bi \left[ 1 + \frac{2}{\pi} h  \right]\Pi(\epsilon_{m,e}) - k_{m,e}^{-1}\left[1+5\frac{h}{\pi}\right]\Pi^2(\epsilon_{m,e}) + \bi\left[  k_{m,e}^{-2}-\frac{1}{12} \right]\Pi^3(\epsilon_{m,e}) + {\cal O}(\epsilon_{m,e}^3\log\epsilon_{m,e}).
\end{align*}
   Consequently, we see that the resonance frequency $kl$, if solving
   (\ref{eq:gov:ke}), must asymptotically behave as (\ref{eq:ke:asy}) for
   $\epsilon_{m,e}\ll 1$. As for the existence of such solutions, one just notices
   that when $kl$ lies in $D_h=\{k\in\mathbb{C}:|kl-k_{m,e}|\leq
   h^{1/2}\}\subset {\cal S}$, then on the boundary of this disk
   \begin{align*}
     &\Bigg|(e^{\bi k l} + 1) - (e^{\bi k l}-1)\left[  c_{00} + <({\rm Id}-{\cal A}_h)^{-1}\{c_{0m}\},\{c_{m0}\}>_{\ell^2}\right] + 2\bi(kl-k_{m,e})\Bigg|\\
     =& {\cal O}(h)\leq 2\sqrt{h}=|-2\bi(kl-k_{m,e})|.
   \end{align*}
   Rouch\'e's theorem indicates that there exists a unique solution to
   (\ref{eq:gov:ke}) in $D_h$.
  \end{proof}
\end{mytheorem}

\subsection{Odd modes}
Suppose now $u$ is odd about $x_3=-l/2$, i.e.,
$u(x_1,x_2,x_3)=-u(x_1,x_2,-l-x_3)$. The original problem can be equivalently
characterized as: find a solution $u\in H^{1}_{\rm loc}(\Omega_h^+)$ solving
\begin{align}
  \label{eq:gov:o}
  \Delta u + k^2 u &= 0,\quad {\rm on}\ \Omega_h^+,\\
  \label{eq:bc:o1}
  \partial_{\nu} u &= 0,\quad {\rm on}\ \partial\Omega_h^+\backslash\overline{\Gamma_l},\\
  \label{eq:bc:o2}
  u &= 0,\quad {\rm on}\ \Gamma_l,
\end{align}
where $\Gamma_l=G_1\times\{x_3=-l/2\}$. As the analysis
follows exactly the same approach as in the even case, we here briefly show the
results.

In $V_h^+$, $u$ could be expressed as
\begin{equation}
  \label{eq:u-:o}
  u(x) = \sum_{m=0}^{+\infty}b_{m}\phi_{m}(x;h)[e^{\bi s_{m}(x_3+l)} -e^{-\bi s_{m} x_3}],
\end{equation}
so that its normal derivative on $\Gamma_h$ becomes
\begin{equation}
  \label{eq:normder:o}
  \partial_{\nu}u(x) = \partial_{x_3}u(x_1,x_2,0) = \sum_{m=0}^{+\infty}b_{m}\bi
  s_{m}\phi_{m}(x;h) [e^{\bi s_{m}l} + 1]\in \tilde{H}^{-1/2}(\Gamma_h).
\end{equation}
In the upper half plane $\mathbb{R}_+^3$, the Neumann data
$\partial_{\nu}u(x)=-\partial_{x_3}u(x)$ in $\tilde{H}^{-1/2}(\Gamma_h)$
uniquely determines the following outgoing solution
\[
  u(x) = -\int_{\Gamma_h} \frac{e^{\bi k |x-y|}}{2\pi|x-y|}\partial_{x_3} u(y)dS(y)\in H^{1}_{\rm loc}(\mathbb{R}_+^3),\quad{\rm on}\quad\mathbb{R}_+^3,
\]
with trace $u|_{\Gamma_h} = -{\cal S}\partial_{x_3}u\in H^{1/2}(\Gamma_h)$ by Lemma~\ref{lem:propSR}.
By the continuity of $\partial_{\nu}u$ on $\Gamma_h$, to ensure that $u\in H^1_{\rm loc}(\Omega_h^+)$, we require that
\[
  u|_{\Gamma_h} = -{\cal S}\partial_{x_3} u =
  \sum_{m=0}^{\infty}b_{m}\phi_{m}(x;h)[e^{\bi s_{m} l} -1]\in H^{1/2}(\Gamma_h),
\]
which is equivalent to for any nonnegative integers $m$,
\[
  (-{\cal S}\partial_{x_3} u,\phi_{m}(x_1,x_2;h))_{L^2(\Gamma_h)} = b_{m}[e^{\bi s_{m}l}-1].
\]
Using again the notation $a_{m}=\lambda_m^{1/4}b_{m}$ for any integer $m>0$, the above equations can be rewritten as the
following infinite number of linear equations in terms of unknowns
$\{b_{0},\{a_{m}\}\}\in\ell^2$, i.e.,
\begin{align}
  \label{eq:b0:o}
  b_{0}[e^{\bi s_{0}l}-1] =&- b_{0}\bi s_{0}h[e^{\bi s_{0}l} + 1]d_{00} - \sum_{m'>0}a_{m'}\lambda_{m'}^{-1/4}\bi s_{m'}h[e^{\bi s_{m'}l} + 1]d_{m'0},\\
  \label{eq:am:o}
  \lambda_{m}^{-1/4}a_{m}[e^{\bi s_{m}l}-1] =&- b_{0}\bi s_{0}h[e^{\bi s_{0}l} + 1]d_{0m}- \sum_{m'>0}a_{m'}\lambda_{m'}^{-1/4}\bi s_{m'}h[e^{\bi s_{m'}l} + 1]d_{m'm},
\end{align}
where $\{d_{m'm}\}_{m,m'=0}^{\infty}$ have been defined in
(\ref{eq:dm'm}). Now set for $m,m'>0$ that
\begin{align}
  c_{00}^{o} =& -\bi \epsilon d_{00},\\
  c_{m'0}^o =& -\lambda_{m'}^{-1/4}\bi s_{m'}h[e^{\bi s_{m'}l} + 1] d_{m'0},\\
  c_{0m}^o =&- \lambda_{m}^{1/4}\frac{\bi \epsilon}{e^{\bi s_{m}l}-1}d_{0m},\\
  c_{m'm}^o =&-\lambda_{m'}^{-1/4}\lambda_{m}^{1/4}\bi s_{m'}h\frac{e^{\bi s_{m'}l} + 1}{e^{\bi s_{m}l}-1}d_{m'm},
\end{align}
we could rewrite the previous linear equations (\ref{eq:b0}) and (\ref{eq:am})
more compactly in the following form,
\begin{align}
  \label{eq:cpb:o}
  b_{0}[e^{\bi s_{0}l}-1] =&b_{0}[e^{\bi s_{0}l} + 1]c^o_{00} + <\{a_{m'}\},\{c^o_{m'0}\}>_{\ell^2},\\
  \label{eq:cpam:o}
  \{a_{m}\}=&b_{0}[e^{\bi s_{0}l} + 1]\{c^o_{0m}\}+ {\cal A}_h^o\{a_{m}\},
\end{align}
where the operator ${\cal A}_h^o$ is defined as: for any $\{f_{m}\}\in\ell^2$,
\begin{align*}
  {\cal A}_h^o\{f_{m}\} = \{\sum_{m>0} c^o_{m'm}f_{m'}\}_{m>0},
\end{align*}
Similar to Lemma~\ref{lem:cm'mA}, we have the following properties.
\begin{mylemma}
  \label{lem:cm'mA:o}
  For $h\ll 1$ and $k\in{\cal B}$:
  \begin{itemize}
  \item[1.]
      \[
        c^o_{00} = -({\cal S}_01,1)_{L^2(\Gamma_1)}\epsilon\bi+ \frac{\epsilon^2}{2\pi} + {\cal O}(\epsilon^3);
      \]
    \item[2.] When $m>0$,
      \[
        c^o_{0m}= \bi\epsilon({\cal S}_{0}1,\lambda_m^{1/4}\phi_{m})_{L^2(\Gamma_1)} +\bi\epsilon^3({\cal R}_01,\lambda_{m}^{1/4}\phi_{m})_{L^2(\Gamma_1)},
      \]
      and $\{c^o_{0m}\}_{m>0}\in \ell^2$;
    \item[3.] When $m'>0$,
      \begin{align*}
        c^o_{m'0}=& \pi({\cal S}_{0} 1,\lambda_{m'}^{1/4}\phi_{m'})_{L^2(\Gamma_1)}+\epsilon^2({\cal R}_0 1, \lambda_{m'}^{1/4}\phi_{m'})_{L^2(\Gamma_1)} + ({\cal S}_0 1, \lambda_{m'}^{-3/4}\phi_{m'})_{L^2(\Gamma_1)}{\cal O}(\epsilon^2),
      \end{align*}
      and $\{c_{m'0}\}_{m'>0}\in \ell^2$;
    \item[4.] When $m',m>0$,
      \begin{align*}
        c^o_{m'm} =& -({\cal S}_0\lambda_{m'}^{1/4}\phi_{m'},\lambda_{m}^{1/4}\phi_{m})_{L^2(\Gamma_1)}-\epsilon^2({\cal R}_0\lambda_{m'}^{1/4}\phi_{m'},\lambda_m^{1/4}\phi_{m})_{L^2(\Gamma_1)} +({\cal S}_0\lambda_{m'}^{-3/4}\phi_{m'},\lambda_m^{1/4}\phi_{m})_{L^2(\Gamma_1)}{\cal O}(\epsilon^2),
      \end{align*}
      and the operator ${\cal A}_h^o$ defined by
      $\{c_{m'm}\}_{m,m'>0}$ is bounded from $\ell^2$ to
      $\ell^2$ and can be decomposed as
      \begin{align}
        {\cal A}_h^o = {\cal P} + \epsilon^2{\cal Q}_h^o,
      \end{align}
      where ${\cal Q}_h^o = \epsilon^{-2}({\cal A}_h^o - {\cal P})$ is uniformly
      bounded from $\ell^2$ to $\ell^2$ for all $k\in{\cal B}$ and $h\ll 1$.
  \end{itemize}
\end{mylemma}
Since ${\cal I}-{\cal P}$ has a bounded inverse, we could obtain the following
theorem, by analogy to Theorem~\ref{thm:PAbi}.
\begin{mytheorem}
  \label{thm:PAbi:o}
  For any $k\in{\cal B}$ and $h\ll1$, ${\cal I}-{\cal A}_h^o$ has a bounded inverse. In fact,
  \[
    ||({\cal I}-{\cal A}_h^o)^{-1} - ({\cal I}-{\cal P})^{-1}|| = {\cal
      O}(\epsilon^2),\quad \epsilon\ll1.
  \]
\end{mytheorem}
\subsubsection{Resonance Frequencies}
Based on Theorem~\ref{thm:PAbi:o}, (\ref{eq:cpb:o}) and (\ref{eq:cpam:o}) are reduced
to the following single equation for the unknown $b_{0}$.
\begin{align}
  \left[  (e^{\bi k l}-1)-(e^{\bi k l}+1)\left(c_{00}^o + <({\cal I}-{\cal A}_h)^{-1}\{c^o_{0m}\},\{c^o_{m0}\}>_{\ell^2}  \right)\right]b_{0} = 0.
\end{align}
Now, based on Lemma~\ref{lem:cm'mA:o} and Theorem~\ref{thm:PAbi:o}, we obtain
the following theorem.
\begin{mytheorem}
  \label{thm:oddres}
  For any width $h\ll 1$, the governing equations (\ref{eq:gov:o}-\ref{eq:bc:o2}) possess nonzero solutions for $k\in{\cal B}$, if and only if
  the following nonlinear equation of $k$
  \begin{equation}
    \label{eq:gov:ko}
    (e^{\bi k l}-1)-(e^{\bi k l}+1)\left(c_{00} + <({\cal I}-{\cal A}_h^o)^{-1}\{c_{0m}^o\},\{c_{m0}^o\}>_{\ell^2} \right)= 0,
  \end{equation}
  has solutions in ${\cal B}$. In fact, these solutions (the so-called resonance
  frequencies) are
  \begin{equation}
    \label{eq:ko:asy}
    kl = k_{m,o} -2\bi \Pi(\epsilon_{m,o}) + 2k_{m,o}^{-1}\Pi(\epsilon_{m,o})(\pi^{-1}\epsilon_{m,o} + 2\Pi(\epsilon_{m,o}) ) + {\cal O}(\epsilon_{m,o}^3), m=1,2,\cdots.
  \end{equation}
  where $k_{m,o} = 2m\pi$ is a Fabry-P\'erot frequency and
  $\epsilon_{m,o}=k_{m,o}h\ll 1$.
  \begin{proof}
    For $h\ll 1$, $\epsilon\ll1$ so that by Theorem~\ref{thm:PAbi} and
    Lemma~\ref{lem:cm'mA}, equation (\ref{eq:gov:ke}) can be reduced to
    \[
      e^{\bi kl}-1 = (e^{\bi kl}+1)\Pi(\epsilon) + {\cal O}(\epsilon^3),
    \]
    which is equivalent to
    \[
      e^{\bi kl} -1 = \frac{2\Pi(\epsilon)}{1-\Pi(\epsilon)} + {\cal O}(\epsilon^3).
    \]
    As the right-hand side approaches $0$ as $\epsilon\to 0$, we see that the
 resonance frequencies must satisfy: for some $m=1,\cdots$, $\delta_{m,o}:= k l - k_{m,o} = o(1)$. Thus,
 \begin{align*}
   \epsilon - \epsilon_{m,o} =& h\delta_{m,o},%\\
   %\epsilon\log\epsilon - \epsilon_{m,e}\log\epsilon_{m,e}  =& h\delta_{m,e} \log\epsilon + h^2\delta_{m,e} + {\cal O}(h^3\delta_{m,e}^2k_{m,e}^{-1}),
 \end{align*}
 as $h \to 0^+$. Note that we cannot allow $m=0$ since $0\notin{\cal B}$. The
 proof follows from the same arguments as in the proof of
 Theorem~\ref{thm:evenres}.
\end{proof}
\end{mytheorem}
\subsection{Two examples and Quality factor}
To conclude this section, we consider two particular shapes for $G_1$, and
shall make a conclusion about the so-called quality factor $Q=-\frac{{\rm
    Re}(k)}{2{\rm Im}(k)}$ for some resonance frequency $k$ in the form of
(\ref{eq:ke:asy}) or (\ref{eq:ko:asy}).

\noindent{\bf Example 1.} When $V_{1,h}$ is generated by a unit square $G_1=(-1/2,1/2)\times(-1/2,1/2)$, we
  can choose for any integer $n\geq 0$,
Let \begin{align}
  \phi_n(x_1;h) =& \left\{
  \begin{array}{lc}
    \sqrt{\frac{1}{h}}&n=0;\\
    \sqrt{\frac{2}{h}}\cos\frac{n\pi x_1}{h}&0<n\mid 2;\\
    \sqrt{\frac{2}{h}}\sin\frac{n\pi x_1}{h}&n\nmid 2,
  \end{array}
  \right.
\end{align}
and $\{\phi_{mn}(x;h) = \phi_m(x_1;h)\phi_n(x_2;h)\}_{m,n=0}^{\infty}$ form a complete and
orthonormal basis in $L^2(\Gamma_h)$. In this case, following \cite{zhlu21}, we in
fact can get by the method of Fourier transform the following identity,
\[
  ({\cal S}_0 1, 1)_{L^2(\Gamma_1)} = \frac{2}{\pi}\log(\sqrt{2}+1)-\frac{2}{3\pi}(\sqrt{2}-1),
\]
to get rid of one unknown constant in (\ref{eq:Pi}).

\noindent{\bf Example 2.} When $V_{1,h}$ is generated by a disk
$G_1=\{x\in\mathbb{R}^2:|x|<1/\sqrt{\pi}\}$ of area $1$,
we can choose for any integer $m, n\geq 0$,
\begin{align}
  \phi_{mne}(r\cos\theta,r\sin\theta;h) =& J_n(h^{-1}\alpha_{mn}r\sqrt{\pi})\cos(n\theta),\\
  \phi_{mno}(r\cos\theta,r\sin\theta;h) =& J_n(h^{-1}\alpha_{mn}r\sqrt{\pi})\sin(n\theta),
\end{align}
where $J_n$ is the $n$-th Bessel function of the first kind, $\alpha_{mn}$ is
the $m$-th smallest root of the equation $J_n'(r)=0$.
$\{\phi_{mne},\phi_{mno}\}_{m,n=0}^{\infty}$ form a complete and orthonormal
basis in $L^2(\Gamma_h)$. This case has been studied by \cite{liazou20}.
%Again, following \cite{}, we in fact can get by the method of Fourier transform the
%following identity,
%\[
%  ({\cal S}_0 1, 1) = \frac{2}{\pi}\log(\sqrt{2}+1)-\frac{2}{3\pi}(\sqrt{2}-1).
%\]

Whatever the shape of $G_1$ is, we in fact can conclude from
Theorems~\ref{thm:evenres} and \ref{thm:oddres} the following result.
\begin{mytheorem}
  For a slab with a single, cylindrical hole $V_{1,h}$ generated by any
  two-dimensional simply-connected Lipschitz domain $G_1$, the quality
  factor $Q$ for the resonance frequency near $m\pi,\ m=1,2,\cdots, $ asymptotically behaves
  as
  \[
    Q=\frac{1}{2h^2m}+{\cal O}(m^{-1}h^{-1}),
  \]
  for $mh\ll 1$. In other words, the leading term of quality factor $Q$ in fact
  is independent of the shape of the cylinder $V_{1,h}$.
\end{mytheorem}

\subsection{Field enhancement}
Suppose now an incident field of a real frequency $k_0$ is specified. If $k_0$
coincides with the real part of some resonance frequency given by
(\ref{eq:ke:asy}) and (\ref{eq:ko:asy}), it is known that the field can be enhanced
inside the slit. Such an anomaly can be simply explained by the proposed
approach. Take the normal incident field $u^{\rm inc}=e^{-\bi k_0 x_3}$ as an
example. The scattering problem can be reduced to two subproblems: (i) with
$u^{\rm inc}/2$ specified in $\mathbb{R}_+^3$, solve (\ref{eq:gov:e}) and
(\ref{eq:bc:e}) for the even field $u^e$; (ii) with $u^{\rm inc}/2$
specified, solve (\ref{eq:gov:o}), (\ref{eq:bc:o1}) and (\ref{eq:bc:o2})
for the odd field $u^o$. The solution to the original problem turns out to
be $u=u^e+u^o$. We consider problem (ii) in the following; problem (i) can be
analyzed similarly. For simplicity, we suppress the superscript $o$. In
$\mathbb{R}_+^3$, define
\begin{align}
  u^{\rm ref}(x) :=u^{\rm inc}(x_1,x_2,x_3)/2 + u^{\rm inc}(x_1,x_2,-x_3)/2=\cos(k_0x_3)/2.
\end{align}
Then, $u-u^{\rm ref}$ is outgoing. Following the same procedures in subsection
2.2, we obtain the following inhomogeneous equation
\begin{align}
-{\cal S}\partial_{x_3} u + u^{\rm ref}(x_1,x_2,0) = \sum_{m=0}^{\infty}b_{m}\phi_{m}(x;h)[e^{\bi s_{m} l} -1],\quad (x_1,x_2)\in \Gamma_{h},
  %&\phi_0(x_1) b_0[e^{\bi s_0 l} - 1]+\sum_{n=1}^{+\infty} n^{-1/2}\phi_n(x_1) a_n[e^{\bi s_n l} - 1]\nonumber\\
  %=&b_0[e^{\bi s_{0} l} + 1]\psi_0(x_1) + \sum_{m=1}^{+\infty}m^{-1/2}a_m[e^{\bi s_{m} l} + 1]\psi_m(x_1) + u^{\rm ref}(x_1,0),\quad |x_1|\leq h/2,
\end{align}
where all definitions remain the same except that $k$ is replaced by $k_0$.
Thus, taking inner product with $\phi_n,n=0,1,\cdots $ yields
\begin{align}
  \label{eq:op1:ince}
  b_{0}[e^{\bi s_{0}l}-1] =&b_{0}[e^{\bi s_{0}l} + 1]c^o_{00} + <\{a_{m'}\},\{c^o_{m'0}\}>_{\ell^2}+b^{\rm ref}_0,\\
  \label{eq:op2:ince}
  \{a_{m}\}=&b_{0}[e^{\bi s_{0}l} + 1]\{c^o_{0m}\}+ {\cal A}_h^o\{a_{m}\}+\{a^{\rm ref}_m\}_{m>0},
\end{align}
where
\begin{align}
  b^{\rm ref}_0 =& \int_{\Gamma_h}u^{\rm ref}(x_1,x_2,0)\phi_0(x;h)dS(x)=h,\\
  a^{\rm ref}_m =& \frac{\lambda_m^{1/4}}{e^{\bi s_m l} -1}\int_{\Gamma_h}u^{\rm ref}(x)\overline{\phi_m(x;h)}dS(x)=0.
\end{align}
For $h\ll 1$, Theorem~\ref{thm:PAbi:o} implies that system
(\ref{eq:op1:ince}-\ref{eq:op2:ince}) can be solved by
\begin{align}
  \label{eq:b0:approx}
    b_0 =& \frac{h}{\left[ (e^{\bi k_0 l} - 1) - (e^{\bi k_0 l}+1)\Pi(k_0h)+{\cal
        O}(k_0^3h^3)\right]} ,\\
  \{a_{m}\}=&b_0(e^{\bi s_0 l} + 1)({\cal I}-{\cal A}_h^{(o)})^{-1}\{c_{0m}^{(o)}\},
\end{align}
For $k_0={\rm Re}(k)$ with $k$ taken as (\ref{eq:ko:asy}) for some
$m\in\mathbb{Z}^+$, $k-k_0={\rm Im}(k)={\cal O}(h^2)$ so
that
\begin{align*}
(e^{\bi k_0 l} -1) - (e^{\bi k_0 l}+1)\Pi(k_0h)
  =& (e^{\bi k l}-1) - (e^{\bi k l}+1)\Pi(k h) + {\cal O}(h^2) ={\cal O}(h^2).
\end{align*}
Thus, $b_0 = {\cal O}(h^{-1})$ and $||\{a_{m}\}||_{\ell^2}={\cal O}(1)$.
We remark that one could follow the proof of Theorem~\ref{thm:evenres} to obtain
the asymptotic behavior of $b_0$ accurate up to ${\cal O}(h^{2})$ as
$h\to 0$. Inside the hole $V_h$, i.e., $x=(x_1,x_2,x_3)\in h\Omega_1\times(-l/2,0)$,
\begin{align*}
  u(x) %=& \sum_{n=0}^{\infty}b_n\phi_n(x_1)\left[ e^{\bi s_n(x_2+l)}-e^{-\bi s_n x_2} \right]\nonumber\\
  =&\frac{1}{h}b_0\left[ e^{\bi k_0(x_3+l)}-e^{-\bi k_0 x_3} \right]+\sum_{m=1}^{\infty}a_m\lambda_{m}^{-1/4}\phi_m(x;h)\left[ e^{\bi s_n(x_3+l)}-e^{-\bi s_m x_3} \right]\nonumber\\
  =&\frac{2}{h}\bi b_0\sin(k_0x_3)+{\cal O}(h^{-1})-\sum_{m=1}^{\infty}a_m\lambda_m^{-1/4}\phi_m(x_1;h)e^{-\bi s_n x_3} + {\cal O}(e^{-\lambda_1 l/(2h)}),
\end{align*}
where we have used the fact that $e^{\bi k_0 l}-1 = {\cal O}(h)$.
Consequently, when $h^2\ll x_3\ll h$, $u(x_1,x_2,x_3)\eqsim {\cal O}(h^{-1})$ and
when $x_3\in(-l/2,0)$ is fixed, $u(x_1,x_2,x_3)\eqsim {\cal O}(h^{-2})$, inducing
field enhancement near the aperture and inside the slit.

\section{Multiple Cylindrical Holes}
In this section, we study resonance frequencies when the slab contains $N$
cylindrical holes $\{V_{i,h}\}_{i=1}^N$. As in \cite{zhlu21}, we begin
with two holes to clarify the main idea.

\subsection{Two holes}
Suppose the slab contains two cylindrical holes $V_{1,h}$ and $V_{2,h}$ centered
at $C_1$ and $C_2$, respectively. Due to the similarity of
even modes and odd modes as discussed before, we here consider the even modes
only and shall directly show the results for odd modes. Let
\begin{align*}
  V_{j,h}^+ =& V_{j,h}\cap\{x\in\mathbb{R}^3:x_3\in(-l/2,0)\},\\
  \Gamma_{j,h} =& (D_j+h G_j)\times\{x_3=0\},\\
  \Gamma_{\cdot,h} =& \cup_{j=1}^2\Gamma_{j,h},\\
  V_{\cdot,h}^+ =& \cup_{j=1}^2V_{j,h}^+,\\
  \Omega_{2,h}^+ =& \mathbb{R}_+^3\cup \Gamma_{\cdot,h}\cup V_{\cdot,h}^+,\\
  \Gamma_{j} = & G_j\times\{x_3=0\},
\end{align*}
and so we need to find $k\in{\cal B}$ such that there exists a nonzero $u\in
H^{1}_{\rm loc}(\Omega_{2,h}^+)$ solving
\begin{align}
  \Delta u + k^2 u = 0,\quad{\rm on}\quad \Omega_{2,h}^+,\\
  \partial_{\nu} u = 0,\quad{\rm on}\quad \partial\Omega_{2,h}^+.
\end{align}
In $V_{j,h}^+,j=1,2$, $u$ can be expressed as
\begin{align}
  \label{eq:u:j}
  u(x) =& \sum_{m=0}^{+\infty}b_{m,j}\phi_{m,j}(x-C_j;h)[e^{\bi s_{m,j}(x_3+l)} + e^{-\bi s_{m,j} x_3}],\quad x\in V_{j,h}^+,\quad j=1,2,
\end{align}
where
\[
  \{\phi_{m,j}(x;h)=h^{-1}\phi_{m,j}((x-C_j)/h)\}_{m=0}^{\infty}
\] are the complete basis in $L^2(\Gamma_{j,h})$,
\[
  s_{m,j} = \sqrt{k^2-\frac{\lambda_{m,j}}{h}},\quad m=0,1,\cdots,
\]
and we recall that $\{\lambda_{m,j}\}$ are the
associated eigenvalues
so that
\begin{align}
  \label{eq:nu:j}
  \partial_{x_3}u(x)=\partial_{\nu}u(x) = \left\{
  \begin{array}{lc}
  \sum_{m=0}^{+\infty}b_{m,1}\bi s_{m,1}\phi_{m}(x-C_1;h) [e^{\bi s_{m,1}l} - 1],& x\in \Gamma_{1,h},\\
  \sum_{m=0}^{+\infty}b_{m,2}\bi s_{m,2}\phi_{m}(x-C_2;h) [e^{\bi s_{m,2}l} - 1],& x\in \Gamma_{2,h},
  \end{array}
\right.
\end{align}
is in $\tilde{H}^{-1/2}(\Gamma_{\cdot,h})$. Now, in the upper half plane
$\mathbb{R}_+^3$, the Neumann data $\partial_{\nu}u(x)=-\partial_{x_3}u(x)$ in
$\tilde{H}^{-1/2}(\Gamma_{\cdot,h})$ uniquely determines the following outgoing solution
\[
  u(x) = -\sum_{j=1}^{2}\int_{\Gamma_{j,h}}\frac{e^{\bi k |x-y|}}{2\pi|x-y|}\partial_{x_3} u(y)dS(y)\in H^{1}_{\rm loc}(\mathbb{R}_+^3),\quad{\rm on}\quad\mathbb{R}_+^3,
\]
Now let
\begin{align}
  \label{eq:int:Sj}
  {\cal S}_j\phi=\int_{\Gamma_{j,h}} \frac{e^{\bi k |x-y|}}{2\pi|x-y|}\phi(y)dS(y),
\end{align}
which is bounded from $\tilde{H}^{-1/2}(\Gamma_{j;h})$ to
$H^{1/2}(\Gamma_{\cdot;h})$ for $j=1,2$. We see on $\Gamma_{\cdot;h}$,
\[
  u = -({\cal S}_1\partial_{\nu}u+{\cal S}_2\partial_{\nu} u),
\]
so that the continuity of $u$ on $\Gamma_{\cdot,h}$
implies that
\begin{align}
  -({\cal S}_1\partial_{\nu} u +{\cal S}_2 \partial_\nu u)|_{\Gamma_{j,h}} =& \sum_{m=0}^{+\infty}b_{m,j}\phi_{m,j}(x-C_j;h)[e^{\bi s_{m,j}l} + 1],\quad{\rm on}\quad\Gamma_{j,h},\quad j=1,2,
  %-({\cal S}_1\partial_{\nu} u +{\cal S}_2 \partial_\nu u)|_{\Gamma_{2,h}} =& \sum_{m=0}^{+\infty}b_{1,m}\phi_{m}(x-C_2;h)[e^{\bi s_{m}l} + 1],\quad{\rm on}\quad\Gamma_{2,h},
  %-({\cal S}_1\partial_{\nu} u +{\cal S}_2 \partial_\nu u)|_{\Gamma_{2,h}} =& \sum_{m, n=0}^{+\infty}b_{2,mn}\phi_{mn}(x-C;h)[e^{\bi s_{mn}l} + 1],\quad{\rm on}\quad\Gamma_{2,h},
\end{align}
which is equivalent to for any integer $m'\geq 0$ that
\begin{align}
  (-({\cal S}_1\partial_{\nu} u +{\cal S}_2 \partial_\nu u)|_{\Gamma_{j,h}},\phi_{m',j}(\cdot-C_j;h))_{L^2(\Gamma_{j,h})} =& b_{m',j}[e^{\bi s_{m',j}l}+1],\quad j=1,2.
  %<-({\cal S}_1\partial_{\nu} u +{\cal S}_2 \partial_\nu u)|_{\Gamma_{2,h}},\phi_{m'n'}(\cdot-C_2;h)> =& b_{2,m'n'}[e^{\bi s_{m'n'}l}+1].
\end{align}
Let
\begin{align}
  \label{eq:dm'm:ij}
d_{m'm}^{ij}=&h^{-1}({\cal S}_i\phi_{m',i}(\cdot-C_i;h)|_{\Gamma_{j,h}},\phi_{m,j}(\cdot-C_j;h))_{L^2(\Gamma_{j,h})},i,j=1,2,
  %\label{eq:dm'n'mn:12}
%d_{m'n',mn}^{12}=&h^{-1}({\cal S}_1\phi_{m'n'}(\cdot;h),\phi_{mn}(\cdot-C;h))_{L^2(\Gamma_{2,h})},
\end{align}
and $a_{m,j}=\lambda_{m,j}^{1/4}b_{m,j}$ so that $\{a_{m,j}\}_{m>0}\in
\ell^2$ for $j=1,2$. By analogy to the equations (\ref{eq:cpb}) and
(\ref{eq:cpam}) in the previous section, we could rewrite the above equations
in terms of matrix operators as following
\begin{align}
  \label{eq:b0:2}
  (e^{\bi kl}+1)\left[
  \begin{array}{cc}
    b_{0,1}\\
    b_{0,2}
  \end{array}
 \right] =& (e^{\bi kl}-1)\left[
  \begin{array}{cc}
    c_{00}^{11} & c_{00,00}^{21}\\
    c_{00}^{12} & c_{00}^{22}\\
  \end{array}
 \right]\left[
  \begin{array}{cc}
    b_{0,1}\\
    b_{0,2}
  \end{array}
 \right] + \left[
  \begin{array}{cc}
    \{c_{m'0}^{11}\} & \{c_{m'0}^{21}\}\\
    \{c_{m'0}^{12}\} & \{c_{m'0}^{22}\}\\
  \end{array}
 \right]\left[
  \begin{array}{cc}
    \{a_{m',1}\}\\
    \{a_{m',2}\}
  \end{array}
 \right],\\
  \label{eq:am:2}
  \left[
  \begin{array}{cc}
    {\cal I} & \\
    & {\cal I}\\
  \end{array}
 \right]\left[
  \begin{array}{cc}
    a_{m,1}\\
    a_{m,2}
  \end{array}
 \right] =& (e^{\bi kl}-1)\left[
  \begin{array}{cc}
    \{c_{0m}^{11}\} & \{c_{0m}^{21}\}\\
    \{c_{0m}^{12}\} & \{c_{0m}^{22}\}\\
  \end{array}
 \right]\left[
  \begin{array}{cc}
    b_{0,1}\\
    b_{0,2}
  \end{array}
 \right] + \left[
  \begin{array}{cc}
    {\cal A}_h^{11} & {\cal A}_h^{21}\\
    {\cal A}_h^{12} & {\cal A}_h^{22}\\
  \end{array}
 \right]\left[
  \begin{array}{cc}
    \{a_{m',1}\}\\
    \{a_{m',2}\}
  \end{array}
 \right],
\end{align}
where for $i,j\in\{1,2\}$, $m,m'>0$ that
\begin{align}
  c_{00}^{ij} =& -\bi \epsilon d_{00}^{ij},\\
  c_{m'0}^{ij} =& -\lambda_{m',i}^{-1/4}\bi s_{m',i}h[e^{\bi s_{m',i}l} - 1] d_{m'0}^{ij},\\
  c_{0m}^{ij} =&- \lambda_{m,i}^{1/4}\frac{\bi \epsilon}{e^{\bi s_{m,j}l}+1}d_{0m}^{ij},\\
  c_{m'm}^{ij} =&-\lambda_{m',i}^{-1/4}\lambda_{m,j}^{1/4}\bi s_{m',i}h\frac{e^{\bi s_{m',i}l} - 1}{e^{\bi s_{m,j}l}+1}d_{m'm}^{ij},
\end{align}
and the operators ${\cal A}_{h}^{ij}$ are defined by (\ref{eq:Ah}) with
$c_{m'm}^{ij}$ in place of $c_{m'm}$. Lemma~\ref{lem:dm'm} can be used
to describe the asymptotic behavior of $d_{m'm}^{ij}$ when $i=j$. If $i\neq j$, we have
\begin{mylemma}
  \label{lem:dm'm:ij}
  For $h\ll1$ and $k\in{\cal B}$, for nonnegative integers $m,m'$, when
  $i,j=1,2$ but $i\neq j$, $d_{m'm}^{ij}$ asymptotically behaves as following:
\begin{align}
  \label{eq:dm'm:ij}
  d_{m'm}^{ij} = \left\{
  \begin{array}{ll}
    \frac{e^{\bi k |C_{ ij }|}}{2\pi k|C_{ij}|}\epsilon + {\cal O}(\epsilon^2), & m=m'=0;\\
    \epsilon^2({\cal R}_0^{ij}\phi_{m',i},\phi_{m,j})_{L^2(\Gamma_j)}, & {\rm otherwise},
  \end{array}
  \right.
\end{align}
where the vector $C_{ij}=C_i-C_j$, ${\cal R}_0^{ij}$ is a uniformly bounded operator from
$\tilde{H}^{-1/2}(\Gamma_i)$ to $H^{1/2}(\Gamma_j)$ for $i=1,2$ and $j=3-i$.
\begin{proof}
  Without loss of generality, suppose $i=2$ and $j=1$. According to the
  definition, we have
  \begin{align*}
    d_{m'm}^{21} =& \frac{1}{2\pi}\int_{\Gamma_2}\int_{\Gamma_1}\frac{he^{\bi k|h(x-y)-C_{21}|}}{|h(x-y)-C_{21}|}\phi_{m,1}(y)dS(y)\phi_{m',2}(x)dS(x)\\
    =& \frac{h}{2\pi}\frac{e^{\bi k |C_{21}|}}{|C_{21}|}\int_{\Gamma_2}\int_{\Gamma_1}\phi_{m,1}(y)dS(y)\phi_{m',2}(x)dS(x)\\
    &+\epsilon^2\frac{1}{2\pi}\int_{\Gamma_2}\int_{\Gamma_1}\epsilon^{-1}\left[  \frac{e^{\bi |\epsilon(x-y)-kC_{21}|}}{|\epsilon(x-y)-kC_{21}|}- \frac{e^{\bi k |C_{21}|}}{|C_{21}|}\right]\phi_{m,1}(y)dS(y)\phi_{m',2}(x)dS(x).
  \end{align*}
It is clear that the first integral on the r.h.s is nonzero only when $m=m'=0$,
and the second integral (excluding the prefactor $\epsilon^2$) has a smooth kernel, which and the gradient of which are uniformly bounded in $\Gamma_2\times \Gamma_1$.
\end{proof}
\end{mylemma}
As an immediate consequence, we get the following lemma.
\begin{mylemma}
  \label{lem:cm'mA:2}
  For $h\ll 1$ and $k\in{\cal B}$: when $i=1,2$ and $j=3-i$,
  \begin{itemize}
  \item[1.]
      \[
        c_{00}^{ij} = -\frac{\bi e^{\bi k|C_{ij}|}}{2\pi k|C_{ij}|}\epsilon^2 + {\cal O}(\epsilon^3);
      \]
    \item[2.] When $m>0$,
      \[
        c_{0m}^{ij}= -\bi\epsilon^3({\cal R}_0^{ij}1,\lambda_{m,j}^{1/4}\phi_{m,j})_{L^2(\Gamma_j)},
      \]
      and $\{c_{0m}^{ij}\}_{m>0}\in \ell^2$;
    \item[3.] When $m'>0$,
      \begin{align*}
        c_{m'0}^{ij}=& - \epsilon^2({\cal R}_0^{ij} 1, \lambda_{m',i}^{1/4}\phi_{m',i})_{L^2(\Gamma_i)},
      \end{align*}
      and $\{c_{m'0}^{ij}\}_{m'>0}\in \ell^2$;
    \item[4.] When $m',m>0$,
      \begin{align*}
        c_{m'm}^{ij} =&-\epsilon^2({\cal R}_0^{i}\lambda_{m',i}^{1/4}\phi_{m',i},\lambda_{m,j}^{1/4}\phi_{m,j})_{L^2(\Gamma_j)}.
      \end{align*}
      and the operator ${\cal A}_h^{ij}$ defined by
      $\{c_{m'm}^{ij}\}_{m,m'>0}$ is bounded from $\ell^2$ to
      $\ell^2$ and $||{\cal A}_h^{ij}||={\cal O}(\epsilon^2)$.
  \end{itemize}
  \begin{proof}
    The proof is similar to that of Lemma~\ref{lem:cm'mA}. We omit the
    details here.
  \end{proof}
\end{mylemma}
By Lemma~\ref{lem:cm'mA:2}, equations (\ref{eq:b0:2}) and (\ref{eq:am:2})
can be transformed to the following two-variable equations:
\begin{equation}
  \label{eq:b1200}
  \left\{ (e^{\bi kl}+1){\cal I}_2-(e^{\bi kl}-1){\rm Diag}\{\Pi_1(\epsilon),\Pi_2(\epsilon) \} + (e^{\bi kl}-1)\bi \epsilon^2{\cal M}_2(k) + {\cal E}_2\right\}\left[
    \begin{array}{c}
      b_{0,1}\\
      b_{0,2}
    \end{array}
  \right] = 0,
\end{equation}
where the function $\Pi_j$ can be defined as $\Pi$ in (\ref{eq:Pi}) but with the
hole $V_{1,h}$ replaced by $V_{j,h}$, ${\cal I}_N$ denotes the $N\times N$
matrix for $N=2$ and
\[
  {\cal M}_2(k) = \left[
    \begin{array}{cc}
      0 & \frac{e^{\bi k|C_{ij}|}}{2\pi k|C_{ij}|}\\
      \frac{e^{\bi k|C_{ij}|}}{2\pi k|C_{ij}|}& 0
    \end{array}
  \right],
\]
and the $2\times 2$ matrix ${\cal E}_2$ consists of elements of ${\cal
  O}(\epsilon^3)$. Note that ${\rm Re}(\Pi_i) = \frac{\epsilon^2}{2\pi}$, which
is independent of $i=1,2$. Now, we characterize the resonance frequencies of
even modes as following.
\begin{mytheorem}
  \label{thm:evenres:2}
 For $h \ll 1$, the resonance frequencies of even modes of the two-hole slab are
 \begin{align}
   \label{eq:asym:ke:2}
   kl = k_{m,e}+\bi\lambda_{j}(\tilde{M}_{2,j}(k_{m,e})-\bi \lambda_j^2(\tilde{M}_{2,j}(k_{m,e})) + {\cal O}(\epsilon_{m,e}^3), \quad j=1,2, m=1,2,\cdots,
 \end{align}
  where  $k_{m,e}=(2m-1)\pi$ is a Fabry-P\'erot frequency, $\epsilon_{m,e}=k_{m,e}h\ll 1$,
 \begin{align}
k_{m,e,j}=&k_{m,e} - 2\bi \Pi_{j}(\epsilon_{m,e}),\\
%e^{\bi k_{m,e,j}l}-1
   \tilde{\cal M}_{2,j}(k_{m,e})=&-2{\rm Diag}\{\Pi_1(k_{m,e,j}h),\Pi_2(k_{m,e,j}h) \}-2\Pi_j(\epsilon_{m,e}){\rm Diag}\{\Pi_1(\epsilon_{m,e}),\Pi_2(\epsilon_{m,e}) \}\nonumber\\
   &+2\bi \epsilon_{m,e}^2{\cal M}_2(k_{m,e}),
   %\tilde{\cal M}_{2,j}(k_{m,e})=&-2{\rm Diag}\{\Pi_1(k_{m,e,j}h),\Pi_2(k_{m,e,j}h) \}-2\Pi_j(\epsilon_{m,e}){\rm Diag}\{\Pi_1(k_{m,e,j}h),\Pi_2(k_{m,e,j}h) \}\nonumber\\
   %&+2\bi \epsilon_{m,e}^2{\cal M}_2(k_{m,e}),
     %\Pi_{j,m,e} =& \Pi(\epsilon_{m,e}) -\bi \epsilon^2_{m,e}\lambda_j({\cal M}_2(k_{m,e})),
     %\delta_{1,j,m,e} =& -\bi \Pi_{1,j,m,e}(\epsilon_{m,e}) - k_{m,e}^{-1}\Pi_{1,j,m,e}^2, \\
     %\Pi_{2,j,m,e} =& \Pi(\epsilon_{m,e}) + \epsilon_{m,e}\lambda_j(S_2(k_{m,e}+\delta_{1,j,m,e},D)),
 \end{align}
 and $\lambda_j(\tilde{\cal M}_{2,j})$ indicates the eigenvalue of $\tilde{\cal
   M}_{2,j}$ closer to $-2\Pi_j(\epsilon_{m,e})$ for $j=1,2$.
 \begin{proof}
   Clearly, (\ref{eq:b1200}) has a nonzero solution $[b_{0,1},b_{0,2}]^{T}$ if and
   only if
   \begin{align}
     \label{eq:mat}
(e^{\bi kl}+1){\cal I}_2-(e^{\bi kl}-1){\rm Diag}\{\Pi_1(\epsilon),\Pi_2(\epsilon) \} + (e^{\bi kl}-1)\bi \epsilon^2{\cal M}_2(k) + {\cal E}_2,
   \end{align}
   has a zero eigenvalue or zero determinant. Since $||(e^{\bi kl}-1)\bi \epsilon^2{\cal M}_2 + {\cal E}_2||_2 = {\cal
     O}(\epsilon^2)$, the resonance frequency $k$ must satisfy
   \begin{align}
     \label{eq:tmp:1}
     (e^{\bi k l}+1) - (e^{\bi k l}-1)\Pi_j(\epsilon) = {\cal O}(\epsilon^2),
   \end{align}
   for some $j=1,2$, since otherwise the matrix in (\ref{eq:mat}) becomes
   diagonally dominant, so that $e^{\bi k l} + 1 = {\cal O}(\epsilon)$. Thus, as
   in Theorem~\ref{thm:evenres}, $kl = k_{m,e} + o(1)$, as $h\to 0$, for some
   $m=1,2,\cdots$. Obviously, $\epsilon\eqsim \epsilon_{m,e}$ so that
   \[
     \delta_{m,e} = kl-k_{m,e}\eqsim (-\bi) \left[  e^{\bi (k-k_{m,e})l} -
       1\right] ={\cal O}(\epsilon) = {\cal O}(\epsilon_{m,e}).
   \]
   Thus, (\ref{eq:tmp:1}) implies that
   \[
     kl = k_{m,e} - 2\bi \Pi_{j}(\epsilon_{m,e}) + {\cal O}(\epsilon_{m,e}^2)
     =k_{m,e,j} + {\cal O}(\epsilon_{m,e}^2),
   \]
   so that we enforce
   \[
(e^{\bi kl}+1){\cal I}_2-\tilde{M_2}(k_{m,e},k_{m,e,j}) + \tilde{\cal E}_2,
   \]
   has a zero eigenvalue, where elements of the $2\times 2$ matrix $\tilde{E}_2$
   are ${\cal O}(\epsilon_{m,e}^3)$. Consequently, we must have that
\begin{equation}
  e^{\bi kl}+1 = \lambda_j(\tilde{\cal M}_2(k_{m,e},k_{m,e,j})) + {\cal O}(\epsilon_{m,e}^3),
\end{equation}
where $\lambda_j$ denotes the eigenvalue (in descending order of
 real part) of $\tilde{\cal M}_2$ closer to $2\Pi_j(\epsilon_{m,e})$; obviously,
 $\lambda_j(\tilde{\cal M}_2) = {\cal O}(\epsilon_{m,e})$. Thus, we have
 \begin{equation}
   \delta_{m,e}^2 - \bi \delta_{m,e} -\lambda_j(\tilde{\cal M}_2(k_{m,e},k_{m,e,j})) + {\cal O}(\epsilon_{m,e}^3) = 0,
 \end{equation}
 so that
 \begin{align*}
   \delta_{m,e} =& \frac{\bi\left(1 -\sqrt{1-4\lambda_{j}(\tilde{\cal M}_2(k_{m,e},k_{m,e,j}))}\right)}{2} + {\cal O}(\epsilon_{m,e}^3)\\
   =&\bi\lambda_{j}(\tilde{\cal M}_2(k_{m,e},k_{m,e,j}))-\bi \lambda_j^2(\tilde{\cal M}_2(k_{m,e},k_{m,e,j})) + {\cal O}(\epsilon_{m,e}^3).
 \end{align*}
We now prove the existence of the two solutions. Assume that $k$ lies in the
   disk $D_h=\{k\in\mathbb{C}:|kl-k_{m,e}|\leq h^{1/2}\}\subset {\cal S}$. Then, on the boundary
   of $D_h$, all entries of
   \[
     (e^{\bi kl}-1)\bi \epsilon^2{\cal M}_2 + {\cal E}_2
   \]
   are ${\cal O}(h^2)$,
   so that by the linearity of determinant,
   \[
     \left|{\rm Det}_1 - {\rm Det}_2\right|={\cal O}(h^{5/2})\leq {\cal O}(h)=\left| {\rm Det}_2 \right|,
   \]
   where
   \begin{align*}
     {\rm Det}_1 = &\Bigg|  (e^{\bi kl}+1){\cal I}_2-(e^{\bi kl}-1){\rm Diag}\{\Pi_1(\epsilon),\Pi_2(\epsilon) \} -
\epsilon (e^{\bi kl}-1){\cal M}_2(k,D) - (e^{\bi kl}-1) {\cal E}_2(\epsilon,D)\Bigg|,\\
     {\rm Det}_2 = &\Bigg|  (e^{\bi kl}+1){\cal I}_2-(e^{\bi kl}-1){\rm Diag}\{\Pi_1(\epsilon),\Pi_2(\epsilon) \}\Bigg|.
   \end{align*}
   For either $j=1,2$, it is clear that on the boundary of $D_h$,
   \begin{align*}
     \left|  2(e^{\bi k l}+1) - (e^{\bi k l}-1)\Pi(\epsilon)+ 2\bi (kl - k_{m,e}) \right| ={\cal O}(h)\leq |-2\bi(kl-k_{m,e})|.
   \end{align*}
   The above two inequalities and Rouch\'e's theorem indicate that there
   are exactly two solutions in $D_h$.
 \end{proof}
\end{mytheorem}

The following theorem characterizes resonance frequencies of odd modes, i.e.,
when the field $u$ satisfies $u(x_1,x_2,-x_3+l/2)=-u(x_1,x_2,x_3+l/2)$.
\begin{mytheorem}
  \label{thm:oddres:2}
 For $h \ll 1$, the resonance frequencies of even modes of the two-hole slab are
 \begin{align}
   \label{eq:asym:ke:2}
   kl = k_{m,o}+\bi\lambda_{j}(\tilde{\cal M}_{2,j}(k_{m,o}))-\bi \lambda_j^2(\tilde{\cal M}_{2,j}(k_{m,o})) + {\cal O}(\epsilon_{m,o}^3), \quad j=1,2, m=1,2,\cdots,
 \end{align}
  where  $k_{m,o}=2m\pi$ is a Fabry-P\'erot frequency, $\epsilon_{m,o}=k_{m,o}h\ll 1$,
 \begin{align}
k_{m,o,j}=&k_{m,o} - 2\bi \Pi_{j}(\epsilon_{m,o}),\\
   \tilde{\cal M}_{2,j}(k_{m,o})=&-2{\rm Diag}\{\Pi_1(k_{m,o,j}h),\Pi_2(k_{m,o,j}h) \}-2\Pi_j(\epsilon_{m,o}){\rm Diag}\{\Pi_1(\epsilon_{m,o}),\Pi_2(\epsilon_{m,o}) \}\nonumber\\
   &+2\bi \epsilon_{m,o}^2{\cal M}_2(k_{m,o}),
                                   %(-2-2\Pi_j(\epsilon_{m,o})){\rm Diag}\{\Pi_1(k_{m,o,j}h),\Pi_2(k_{m,o,j}h) \} +2\bi \epsilon_{m,o}^2{\cal M}_2(k_{m,o}),
     %\Pi_{j,m,e} =& \Pi(\epsilon_{m,e}) -\bi \epsilon^2_{m,e}\lambda_j({\cal M}_2(k_{m,e})),
     %\delta_{1,j,m,e} =& -\bi \Pi_{1,j,m,e}(\epsilon_{m,e}) - k_{m,e}^{-1}\Pi_{1,j,m,e}^2, \\
     %\Pi_{2,j,m,e} =& \Pi(\epsilon_{m,e}) + \epsilon_{m,e}\lambda_j(S_2(k_{m,e}+\delta_{1,j,m,e},D)),
 \end{align}
 and $\lambda_j(\tilde{\cal M}_{2,j})$ indicates the eigenvalue of $\tilde{\cal
   M}_{2,j}$ closer to $-2\Pi_j(\epsilon_{m,o})$ for $j=1,2$.
 \begin{proof}
   The proof follows from similar arguments as in Theorem~\ref{thm:evenres:2}.
 \end{proof}
 \end{mytheorem}
 The above results can be readily extended to a slab with the $N$ holes
 $\{V_{j,h}\}_{j=1}^N$ centered at $\{C_j\}_{j=1}^N$. We
 state our main result in the following.
 \begin{mytheorem}
   \label{thm:res:N}
 For $h \ll 1$, the resonance frequencies of a slab containing $\{V_{j,h}\}_{j=1}^N$ are
 \begin{align}
   \label{eq:asym:k:N}
   kl = k_{m}+\bi\lambda_{j}(\tilde{\cal M}_{N,j}(k_{m}))-\bi \lambda_j^2(\tilde{\cal M}_{N,j}(k_{m})) + {\cal O}(\epsilon_{m}^3),
 \end{align}
 for $j=1,\cdots, N$ and $m=1,2,\cdots,$ where $k_{m}=m\pi$ is a Fabry-P\'erot
 frequency, $\epsilon_{m}=k_{m}h\ll 1$,
 \begin{align}
k_{m,j}=&k_{m} - 2\bi \Pi_{j}(\epsilon_{m}),\\
   {\cal M}_N(k,\{C_i\}_{i=1}^{N}) =&\left[
  \begin{array}{cccc}
    0 &  \frac{e^{\bi k|C_{12}|}}{2\pi k|C_{12}|} & \cdots &  \frac{e^{\bi k|C_{1N}|}}{2\pi k|C_{1N}|}\\
    \frac{e^{\bi k|C_{21}|}}{2\pi k|C_{21}|} & 0 & \cdots & \frac{e^{\bi k|C_{2N}|}}{2\pi k|C_{2N}|} \\
    \vdots & \vdots & \ddots & \vdots\\
    %H_0^{(1)}(kD_{N-1,1}) & H_0^{(1)}(kD_{N-1,2})&\cdots & 0 & H_0^{(1)}(kD_{N-1,N})\\
    \frac{e^{\bi k|C_{N1}|}}{2\pi k|C_{N1}|} & \frac{e^{\bi k|C_{N2}|}}{2\pi k|C_{N2}|}&\cdots & 0\\
  \end{array}
\right],\\
   \tilde{\cal M}_{N,j}(k_{m})=&-2{\rm Diag}\{\Pi_1(k_{m,o,j}h),\cdots,\Pi_N(k_{m,o,j}h) \}\\
   &-2\Pi_j(\epsilon_{m,o}){\rm Diag}\{\Pi_1(\epsilon_{m,o}),\cdots,\Pi_N(\epsilon_{m,o}) \}+2\bi \epsilon_{m,o}^2{\cal M}_N(k_{m,o}),
%
%                                 (-2-2\Pi_j(\epsilon_{m})){\rm Diag}\{\Pi_1(k_{m,j}h),\cdots, \Pi_N(k_{m,j}h) \} +2\bi \epsilon_{m}^2{\cal M}_N(k_{m}),
     %\Pi_{j,m,e} =& \Pi(\epsilon_{m,e}) -\bi \epsilon^2_{m,e}\lambda_j({\cal M}_2(k_{m,e})),
     %\delta_{1,j,m,e} =& -\bi \Pi_{1,j,m,e}(\epsilon_{m,e}) - k_{m,e}^{-1}\Pi_{1,j,m,e}^2, \\
     %\Pi_{2,j,m,e} =& \Pi(\epsilon_{m,e}) + \epsilon_{m,e}\lambda_j(S_2(k_{m,e}+\delta_{1,j,m,e},D)),
 \end{align}
 and $\lambda_j(\tilde{\cal M}_{N,j})$ indicates the eigenvalue of $\tilde{\cal M}_{N,j}$ closest to $-2\Pi_j(\epsilon_{m})$ for $j=1,2$.

 \begin{proof}
   The proof is analogous to that of Theorems~\ref{thm:evenres:2} and \ref{thm:oddres:2}.
 \end{proof}
\end{mytheorem}
\begin{myremark}
  When all $\{V_{j,h}\}_{j=1}^N$ are generated by the same Lipschitz domain, say
  $G_1$, we could simplify the above formulae and get:
\begin{align}
   \label{eq:kl:N'}
   kl = &k_{m} -2\bi \Pi_{j,m}(\epsilon_{m}) + 2k_{m}^{-1}\Pi_{j,m}(\epsilon_{m})(\pi^{-1}\epsilon_{m} + 2\Pi_{j,m}(\epsilon_{m}) ) + {\cal O}(\epsilon_{m}^3),\quad j=1,\cdots, N, m=1,2,\cdots,
   %kl = &k_{m}-\bi \Pi_{1,j,m} - k_{m}^{-1}\Pi_{1,j,m}^2 + {\cal O}(\epsilon_{m}^2\log\epsilon_{m}),j=1,\cdots, N, m=1,2,\cdots,
 \end{align}
 where $k_m=m\pi$ is a Fabry-P\'erot frequency, $\epsilon_m=k_mh\ll 1$,
 \begin{align*}
     \Pi_{j,m} =& \Pi_1(\epsilon_{m}) -\bi \epsilon_m^2\lambda_j({\cal M}_N(k_{m},\{C_i\}_{i=1}^{N})),
 \end{align*}
and $\lambda_j({\cal M}_N(k,\{C_i\}_{i=1}^N))$ indicates the $j$-th eigenvalue (in
descending order of real part) of ${\cal M}_N$. Consequently, the quality factor
$Q$ for the resonance frequency $k$ in (\ref{eq:kl:N'}) behaves as
\[
  Q = \frac{1}{(2+2{\rm Im}(\lambda_j)\pi)mh^2} + {\cal O}(m^{-1}h^{-1}),\quad
  mh\ll 1.
\]
Clearly, the leading behavior of $Q$ does not rely on the choice of shape
of $\Gamma_1$, but only the locations of $V_{j,h}$ as $h\ll1$.
\end{myremark}
\begin{myremark}
  In fact, all the previous theoretical results can be directly genearlized to any
  dimensions greater than three.
\end{myremark}

The field enhancement in the $N$-hole slab can be analyzed by similar arguments as in section 2.4. We omit the details here.
\section{Conclusion}
This paper has developed a simple Fourier-matching method to rigorously study
resonance frequencies of a sound-hard slab with a finite number of arbitrarily
shaped cylindrical holes of diameter ${\cal O}(h)$ for $h\ll1$. Outside the
holes, a sound field was expressed in terms of its normal derivatives on the
apertures of holes. Inside each hole, since the vertical variable can be
separated, the field was expressed in terms of a countable set of Fourier basis
functions. Matching the field on each aperture yields a linear system of
countable equations in terms of a countable set of unknown Fourier coefficients.
The linear system was further reduced to a finite-dimensional linear system by
studying the well-posedness of a closely related boundary value problem for each
hole for $h\ll1$, so that only the leading Fourier coefficient of each hole was
preserved in the final finite-dimensional system. The resonance frequencies are
those making the resulting finite-dimensional linear system rank deficient. By
regular asymptotic analysis for $h\ll1$, we obtained a systematic asymptotic
formula for characterizing the resonance frequencies by the 3D subwavelength
structure. The formula revealed an important fact that when all holes are of the
same shape, the $Q$-factor for any resonance frequency asymptotically behaves as
${\cal O}(h^{-2})$ for $h\ll1$ with its prefactor independent of shapes of
holes. This indicates that the shape of subwavelength structures in fact plays
less significant roles in realizing high-Q resonators.

Since the proposed Fourier matching method does not need to analyze the
complicated Green function of each hole nor need to know the shape of each hole,
we expect that the method can be extended to analyze more complicated and
realistic structures. Our future plan is to extend the current method to analyze
resonances of electro-magnetic scattering problems by 3D subwavelength
structures.
\section*{Acknowledgement}
WL would like to thank Prof. Hai Zhang of Hong Kong University of Science and
Technology for some useful discussions.

\bibliographystyle{plain}
\bibliography{wt}

\begin{thebibliography}{10}

\bibitem{braholsch20}
Brand\ ao~R., Holley~J. R., and Schnitzer O.
\newblock Boundary-layer effects on electromagnetic and acoustic extraordinary
  transmission through narrow slits.
\newblock {\em Proc. R. Soc. A.}, page 20200444, 2020.

\bibitem{astlalpal00}
P.~Astilean, S.~Lalanne and M.~Palamaru.
\newblock Light transmission through metallic channels much smaller than the
  wavelength.
\newblock {\em Opt. Commun.}, 175:265--273, 2000.

\bibitem{babbontri10}
J-F. Babadjian, E.~Bonnetier, and F.~Triki.
\newblock Enhancement of electromagnetic fields caused by interacting
  subwavelength cavities.
\newblock {\em Multiscale Model. Simul.}, 8(4):1383--1418, 2010.

\bibitem{bonsta94}
A.~Bonnet-Bendhia and F.~Starling.
\newblock Guided waves by electromagnetic gratings and nonuniqueness examples
  for the diffraction problem.
\newblock {\em Math. Methods Appl. Sci.}, 17:305--338, 1994.

\bibitem{bontri10}
E.~Bonnetier and F.~Triki.
\newblock Asymptotic of the green function for the diffraction by a perfectly
  conducting plane perturbed by a sub-wavelength rectangular cavity.
\newblock {\em Math. Methods Appl. Sci.}, 33:772--798, 2010.

\bibitem{carmahgar11}
S.~Carretero-Palacios, O.~Mahboub, F.~J. Garcia-Vidal, L.~Martin-Moreno, S.~G.
  Rodrigo, C.~Genet, and T.~W. Ebbesen.
\newblock Mechanisms for extraordinary optical transmission through bull’s
  eye structures.
\newblock {\em Opt. Express}, 19:10429--10442, 2011.

\bibitem{cladurjoltor06}
M.~Clausel, M.~Durufle, P.~Joly, and Tordeux S.
\newblock A mathematical analysis of the reso- nance of the finite thin slots.
\newblock {\em Appl. Numer. Math.}, 56:1432--1449, 2006.

\bibitem{ebblezwol98}
T.~W. Ebbesen, H.~J. Lezec, H.~F. Ghaemi, T.~Thio, and P.~A. Wolff.
\newblock Extraordinary optical transmission through sub-wavelength hole
  arrays.
\newblock {\em Nature}, 391:667–669, 1998.

\bibitem{chenet13}
X.~Chen et~al.
\newblock Atomic layer lithography of wafer-scale nanogap arrays for extreme
  confinement of electro-magnetic waves.
\newblock {\em Nat. Commun.}, 4:2361, 2013.

\bibitem{gaoliyua17}
Y.~Gao, P.~Li, and X.~Yuan.
\newblock Electromagnetic field enhancement in a subwavelength rectangular open
  cavity.
\newblock {\em arXiv:1711.06804}, 2017.

\bibitem{garmarebbkui10}
F.~J. Garcia-Vidal, L.~Martin-Moreno, T.~W. Ebbesen, and L.~K. Kuipers.
\newblock Light passing through subwavelength apertures.
\newblock {\em Rev. Mod. Phys.}, 82:729787, 2010.

\bibitem{genebb07}
C.~Genet and T.~W. Ebbesen.
\newblock Light in tiny holes.
\newblock {\em Nature}, 445:39--46, 2007.

\bibitem{holsch19}
J.~R. Holley and Schnitzer O.
\newblock Extraordinary transmission through a narrow slit.
\newblock {\em Wave Motion}, 91:102381, 2019.

\bibitem{huyualu20}
Z.~Hu, Yuan L., and Y.~Y. Lu.
\newblock Resonant field enhancement near bound states in the continuum on
  periodic structures.
\newblock {\em Physical Review A}, 101(4):043825, 2020.

\bibitem{joltor06a}
P.~Joly and S.~Tordeux.
\newblock Asymptotic analysis of an approximate model for time harmonic waves
  in media with thin slots.
\newblock {\em ESAIM Math. Model. Numer. Anal.}, 40:63--97, 2006.

\bibitem{joltor06b}
P.~Joly and S.~Tordeux.
\newblock Matching of asymptotic expansions for wave propagation in media with
  thin slots i: The asymptotic expansion.
\newblock {\em Multiscale Model. Simul.}, 5:304--336, 2006.

\bibitem{joltor08}
P.~Joly and S.~Tordeux.
\newblock Matching of asymptotic expansions for wave propagation in media with
  thin slots ii: the error estimates.
\newblock {\em ESAIM Math. Model. Numer. Anal.}, 42:193--221, 2008.

\bibitem{liazou20}
Y.~Liang and J.~Zou.
\newblock Acoustic scattering and field enhancement through a single aperture.
\newblock {\em arXiv:2011.05887}, 2020.

\bibitem{lienyllud83}
B.~Liedberg, C.~Nylander, and I.~Lundstrom.
\newblock Surface plasmons resonance for gas detection and biosensing.
\newblock {\em Sensors Actuators}, 4(299), 1983.

\bibitem{linshizha20}
J.~Lin, S.~P. Shipman, and H.~Zhang.
\newblock A mathematical theory for fano resonance in a periodic array of
  narrow slits.
\newblock {\em SIAM J. Math. Analy.}, 80(5):2045--2070, 2020.

\bibitem{linzha17}
J.~Lin and H.~Zhang.
\newblock Scattering and field enhancement of a perfect conducting narrow slit.
\newblock {\em SIAM J. Appl. Math.}, 77(3):951--976, 2017.

\bibitem{linzha18a}
J.~Lin and H.~Zhang.
\newblock Scattering by a periodic array of subwavelength slits i: field
  enhancement in the diffraction regime.
\newblock {\em Multiscale Model. Simul.}, 16(2):922--953, 2018.

\bibitem{linzha18b}
J.~Lin and H.~Zhang.
\newblock Scattering by a periodic array of subwavelength slits ii: surface
  bound state, total transmission and field enhancement in homogenization
  regimes.
\newblock {\em Multiscale Model. Simul.}, 16(2):954--990, 2018.

\bibitem{liulal08}
H.~Liu and P.~Lalanne.
\newblock Microscopic theory of the extraordinary optical transmission.
\newblock {\em Nature}, 452(7188):728--731, 2008.

\bibitem{mcl00}
W.~McLean.
\newblock {\em Strongly Elliptic Systems and Boundary Integral Equations}.
\newblock Cambridge University Press, New York, NY, 2000.

\bibitem{sarvig07}
M.~Sarrazin and J.~P. Vigneron.
\newblock Bounded modes to the rescue of optical transmission.
\newblock {\em Europhysics News}, 38:27--31, 2007.

\bibitem{seoetal09}
M.~A. et~al. Seo.
\newblock Terahertz field enhancement by a metallic nano slit operating beyond
  the skin-depth limit.
\newblock {\em Nat. Photonics}, 3:152--156, 2009.

\bibitem{shi10}
S.~P. Shipman.
\newblock {\em Resonant scattering by open periodic waveguides, Chapter 2 in
  Wave Propagation in Periodic Media: Analysis, Numerical Techniques and
  Practical Applications, M. Ehrhardt, ed., E-Book Series PiCP}, volume~1.
\newblock Bentham Science Publishers, 2010.

\bibitem{shivol07}
S.~P. Shipman and D.~Volkov.
\newblock Guided modes in periodic slabs: existence and nonexistence.
\newblock {\em SIAM J. Appl. Math.}, 67:687--713, 2007.

\bibitem{stupodgor10}
B.~Sturman, E.~Podivilov, and M.~Gorkunov.
\newblock Transmission and diffraction properties of a narrow slit in a perfect
  metal.
\newblock {\em Phys. Rev. B}, 82:115419, 2010.

\bibitem{tak01}
Y.~Takakura.
\newblock Optical resonance in a narrow slit in a thick metallic screen.
\newblock {\em Phys. Rev. Lett.}, 99:5601–5603, 2001.

\bibitem{yansam02}
F.~Yang and J.~R. Sambles.
\newblock Resonant transmission of microwaves through a narrow metallic slit.
\newblock {\em Phys. Rev. Lett.}, 89:063901, 2002.

\bibitem{zhlu21}
J.~Zhou and W.~Lu.
\newblock Numerical analysis of resonances by a slab of subwavelength slits by
  fourier transform-based matching method.
\newblock {\em submitted}, 2021.

\end{thebibliography}
\end{document}